\documentclass[a4paper]{article}\def\zibreport{1}

\usepackage{ifthen}
\usepackage{graphicx}
\usepackage{multicol}
\usepackage{footmisc}
\usepackage{mathdots}
\usepackage{caption}
\usepackage{subcaption}
\usepackage{adjustbox}
\usepackage{amsmath,amsfonts,amssymb,mathrsfs}
\usepackage{tabularx,supertabular,booktabs}
\usepackage[textsize=scriptsize, textwidth=3.3cm]{todonotes}
\usepackage{xspace}
\usepackage{lscape}
\usepackage{multirow}
\usepackage{tabularx}
\usepackage{comment}
\usepackage[utf8]{inputenc}
\usepackage{hyperref}
\usepackage{cleveref}
\usepackage{blindtext}
\usepackage{url}

\usepackage{verbatim}
\usepackage{algorithm}
\usepackage{algorithmic}
\Crefname{ALC@unique}{Line}{Lines}

\usepackage{forest}
\tikzset{
    solid node/.style={circle,draw,inner sep=1.2,fill=black},
    hollow node/.style={circle,draw,inner sep=1.2},
    left label/.style={above left,midway},
    right label/.style={above right,midway}
}

\usepackage{tikz}
\usetikzlibrary{matrix,calc,positioning}
\usetikzlibrary{petri}	
\usetikzlibrary{shapes.misc}
\tikzset{cross/.style={cross out, draw=black, minimum size=2*(#1-\pgflinewidth), inner sep=0pt, outer sep=0pt},
	cross/.default={1pt}}
\usepackage{pgfplots}
\pgfplotsset{compat=1.16}

\newcommand{\myorcidlink}[1]{\,\href{https://orcid.org/#1}{\raisebox{-0.45ex}{\includegraphics[width=1.8ex]{orcid}}}}
\newcommand{\myurl}[1]{\textsf{\footnotesize \url{#1}}\xspace}
\newcommand{\matr}[1]{\begin{bmatrix} #1 \end{bmatrix}}
\newcommand{\allcaps}[1]{\protect\scalebox{0.93}{#1}}

\newcommand{\name}[1]{\mbox{#1}\xspace}
\newcommand{\nameCaps}[1]{\name{\allcaps{#1}}}

\newcommand{\mcol}[3]{\multicolumn{#1}{#2}{#3}}
\newcommand{\mrow}[3]{\multirow{#1}{#2}{$#3$}}

\newcommand{\remix}{\allcaps{REMix}\xspace}
\newcommand{\remixmiso}{\allcaps{REMix-MISO}\xspace}
\newcommand{\unseen}{\allcaps{UNSEEN}\xspace}
\newcommand{\beamme}{\allcaps{BEAM-ME}\xspace}

\newcommand\scalemath[2]{\scalebox{#1}{\mbox{\ensuremath{\displaystyle #2}}}}

\newcommand{\OOPS}{\nameCaps{OOPS}}
\newcommand{\OOQP}{\nameCaps{OOQP}}
\newcommand{\PIPSPETRA}{\nameCaps{PIPS-IPM}}
\newcommand{\JUWELS}{\nameCaps{JUWELS}}

\newcommand{\UNSEEN}{\nameCaps{UNSEEN}}

\newcommand{\CPLEX}{\nameCaps{CPLEX}}

\newcommand{\GUROBI}{\name{Gurobi}}
\newcommand{\PIPS}{\nameCaps{PIPS-IPM++}}

\newcommand{\PARDISO}{\nameCaps{PARDISO}}
\newcommand{\MATwentySeven}{\nameCaps{MA27}}
\newcommand{\MAFiftySeven}{\nameCaps{MA57}}
\newcommand{\MUMPS}{\nameCaps{MUMPS}}
\newcommand{\WSMP}{\nameCaps{WSMP}}
\newcommand{\MAEightySix}{\nameCaps{MA86}}
\newcommand{\MKL}{\nameCaps{MKL}}

\newcommand{\COPT}{\nameCaps{COPT}}

\newcommand{\MadNLP}{\nameCaps{MadNLP}}

\newcommand{\SIMPLE}{\nameCaps{{SIMPLE}}}

\newcommand{\LP}{\nameCaps{LP}}
\newcommand{\LPs}{\nameCaps{LPs}}
\newcommand{\MIP}{\nameCaps{MIP}}

\newcommand{\MIPs}{\nameCaps{MIPs}}

\newcommand{\MPI}{\nameCaps{MPI}}
\newcommand{\HPC}{\nameCaps{HPC}}

\newcommand{\OpenMP}{\nameCaps{OpenMP}}
\newcommand{\KKT}{\nameCaps{KKT}}

\newcommand{\IPM}{\nameCaps{IPM}}
\newcommand{\IPMs}{\nameCaps{IPMs}}

\newcommand{\HSCA}{\nameCaps{HSCA}}
\newcommand{\AHLP}{\nameCaps{AHLP}}
\newcommand{\AHLPs}{\nameCaps{AHLPs}}

\newcommand{\D}{\ensuremath{\Delta}}
\newcommand{\var}{\ensuremath{x}\xspace}

\newcommand{\primalvector}{\var}
\newcommand{\primalvectorCap}{\ensuremath{X}\xspace}
\newcommand{\T}{^{T}}
\newcommand{\inv}{^{-1}}
\newcommand{\vecOnes}{\ensuremath{e}\xspace}

\newcommand{\centralityFactor}{\ensuremath{\tau}\xspace}
\newcommand{\rhsvectorEq}{\ensuremath{b}\xspace}

\newcommand{\coefmatrixEq}{\ensuremath{A}\xspace}
\newcommand{\dualsEq}{\ensuremath{y}\xspace}
\newcommand{\coefmatrixIneq}{\ensuremath{C}\xspace}
\newcommand{\dualsIneq}{\ensuremath{z}\xspace}

\newcommand{\dualsLowerBounds}{\ensuremath{\gamma}\xspace}
\newcommand{\dualsIneqCap}{\ensuremath{Z}\xspace}

\newcommand{\dualsLowerBoundsCap}{\ensuremath{\Gamma}\xspace}

\newcommand{\objectiveLinearVector}{\ensuremath{c}\xspace}

\newcommand{\subjectto}{\mbox{subject to}\xspace}

\newcommand{\bigO}{\mathcal{O}}

\newcommand{\N}{\mathbb{N}}
\newcommand{\R}{\mathbb{R}}
\newcommand{\Rtimes}[2]{\R^{{#1}\times{#2}}}

\newcommand{\Rmn}{\Rtimes{m}{n}}

\DeclareMathOperator{\diag}{diag}

\definecolor{seagreen}{rgb}{0.18,0.74,0.56}
\definecolor{darkgreen}{rgb}{0.0,0.45,0.00}
\definecolor{navyblue}{rgb}{0.0,0.0,0.5}
\definecolor{steelblue}{rgb}{0.27,0.51,0.71}
\definecolor{siennabrown}{rgb}{0.63,0.32,0.18}
\definecolor{firebrickred}{rgb}{0.69,0.13,0.13}
\definecolor{gray75}{rgb}{0.75,0.75,0.75}
\definecolor{orange}{rgb}{.843,0.671,0.078}
\definecolor{gold}{rgb}{1.0,0.84,0.0}
\definecolor{scipyellow}{HTML}{FFFFD6}
\definecolor{soplexred}{HTML}{FFD8D8}
\definecolor{zimplgreen}{HTML}{D8FFD8}
\definecolor{ugblue}{HTML}{CFEFFF}
\definecolor{gcgorange}{HTML}{FFDAB9}

\definecolor{mDarkBrown}{HTML}{604c38}
\definecolor{mDarkTeal}{HTML}{23373b}
\definecolor{mLightBrown}{HTML}{EB811B}
\definecolor{mLightblue}{HTML}{14B03D}

\definecolor{c0}{HTML}{000060}
\definecolor{c1}{HTML}{0000FF}
\definecolor{c2}{HTML}{36648B}
\definecolor{c3}{HTML}{4682B4}
\definecolor{c4}{HTML}{5CACEE}
\definecolor{c5}{HTML}{FF0000}
\definecolor{c6}{HTML}{008888}
\definecolor{c7}{HTML}{00DD99}
\definecolor{c8}{HTML}{527B10}
\definecolor{c9}{HTML}{7BC618}
\definecolor{c10}{HTML}{8DD8F8}
\definecolor{background}{HTML}{FFFFFF}

\definecolor{cp1}{HTML}{F2AF29}
\definecolor{cp2}{HTML}{05668D}
\definecolor{cp3}{HTML}{02C39A}
\definecolor{cp4}{HTML}{4F345A}
\definecolor{cp5}{HTML}{F2B5D4}
\definecolor{cp6}{HTML}{DA2C38}

\newcommand{\PIPSBl}{\ensuremath{F}\xspace}
\newcommand{\PIPSDl}{\ensuremath{G}\xspace}
\newcommand{\PIPSBlloc}{\mathcal{F}}

\newcommand{\PIPSBlglob}{\ensuremath{\textbf{F}}}

\newcommand{\ESM}{\nameCaps{ESM}}
\newcommand{\ESMs}{\nameCaps{ESMs}}

\renewcommand\dualsLowerBounds{\dualsIneq}
\renewcommand\dualsLowerBoundsCap{\dualsIneqCap}

\newcommand{\TheTitle}{A Massively Parallel Interior-Point Method for Arrowhead Linear Programs with Local Linking Structure}
\newcommand{\TheAuthors}{N.--C. Kempke, D. Rehfeldt, T. Koch}

\hypersetup{%
    pdftitle={\TheTitle},
    pdfauthor={\TheAuthors}
}

\ifthenelse{\zibreport = 1}{
    \usepackage{./zibtitlepage}
}{
    \title{{\TheTitle}}

    \headers{A Massively Parallel Interior-Point Method}{\TheAuthors}

    \author{
    Nils--Christian Kempke\thanks{Applied Algorithmic Intelligence Methods Department, Zuse Institute Berlin, Takustr.~7, 14195~Berlin, Germany (\email{kempke@zib.de}, \email{rehfeldt@zib.de}, \email{koch@zib.de}).}
    \and
	Daniel Rehfeldt\footnotemark[1] \thanks{IVU Traffic Technologies, Bundesallee 88, 12161 Berlin, Germany}
	\and
  	Thorsten Koch\footnotemark[1]\ \thanks{Chair of Software and Algorithms for Discrete Optimization, Technische Universit{\"a}t Berlin, Stra{\ss}e des 17.~Juni~136, 10623~Berlin.}
    }
}

\begin{document}

\ifthenelse{\zibreport = 1}{
    \ZTPTitle{\TheTitle}
    \title{\TheTitle}

    \ZTPAuthor{
        \ZTPHasOrcid{Nils--Christian Kempke}{0000-0003-4492-9818},
        \ZTPHasOrcid{Daniel Rehfeldt}{0000-0002-2877-074X},
        \ZTPHasOrcid{Thorsten Koch}{0000-0002-1967-0077}
    }
    \author{
        \ZTPHasOrcid{Nils--Christian Kempke}{0000-0003-4492-9818},\and\
        \ZTPHasOrcid{Daniel Rehfeldt}{0000-0002-2877-074X},\and\
        \ZTPHasOrcid{Thorsten Koch}{0000-0002-1967-0077}
    }

    \ZTPInfo{Preprint}
    \ZTPNumber{24-13}
    \ZTPMonth{December}
    \ZTPYear{2024}

    \date{\normalsize March 4, 2026}
}{}

\maketitle

\begin{abstract}
In practice, non-specialized interior-point algorithms often cannot utilize the massively parallel compute resources offered by modern many- and multi-core compute platforms.
However, efficient distributed solution techniques are required, especially for large-scale linear programs.
This article describes a new decomposition technique for systems of linear equations, implemented in the parallel interior-point solver \PIPS.
The algorithm exploits a matrix structure commonly found in optimization problems: a doubly bordered block-diagonal or arrowhead structure with linking constraints and variables often only linking few, consecutive blocks.
This structure is preserved in the linear KKT systems solved during each iteration of the interior-point method.
We present a hierarchical Schur complement decomposition that distributes and solves the linear optimization problem.
It is designed for high-performance architectures and scales well with the availability of additional computing resources.
The decomposition approach uses the border constraints' locality to decouple the factorization process.
Our approach is motivated by large-scale economic dispatch problems but can also be applied to other problem classes.
We demonstrate the performance of our method on a set of mid- to large-scale instances, some of which have more than $10^{9}$ nonzeros in their constraint matrices.

\ifthenelse{\zibreport = 0}{

}{}

\end{abstract}

\section{Introduction}
\label{sec:intro}

\begin{figure}
    \centering
    \includegraphics[width=0.4\linewidth]{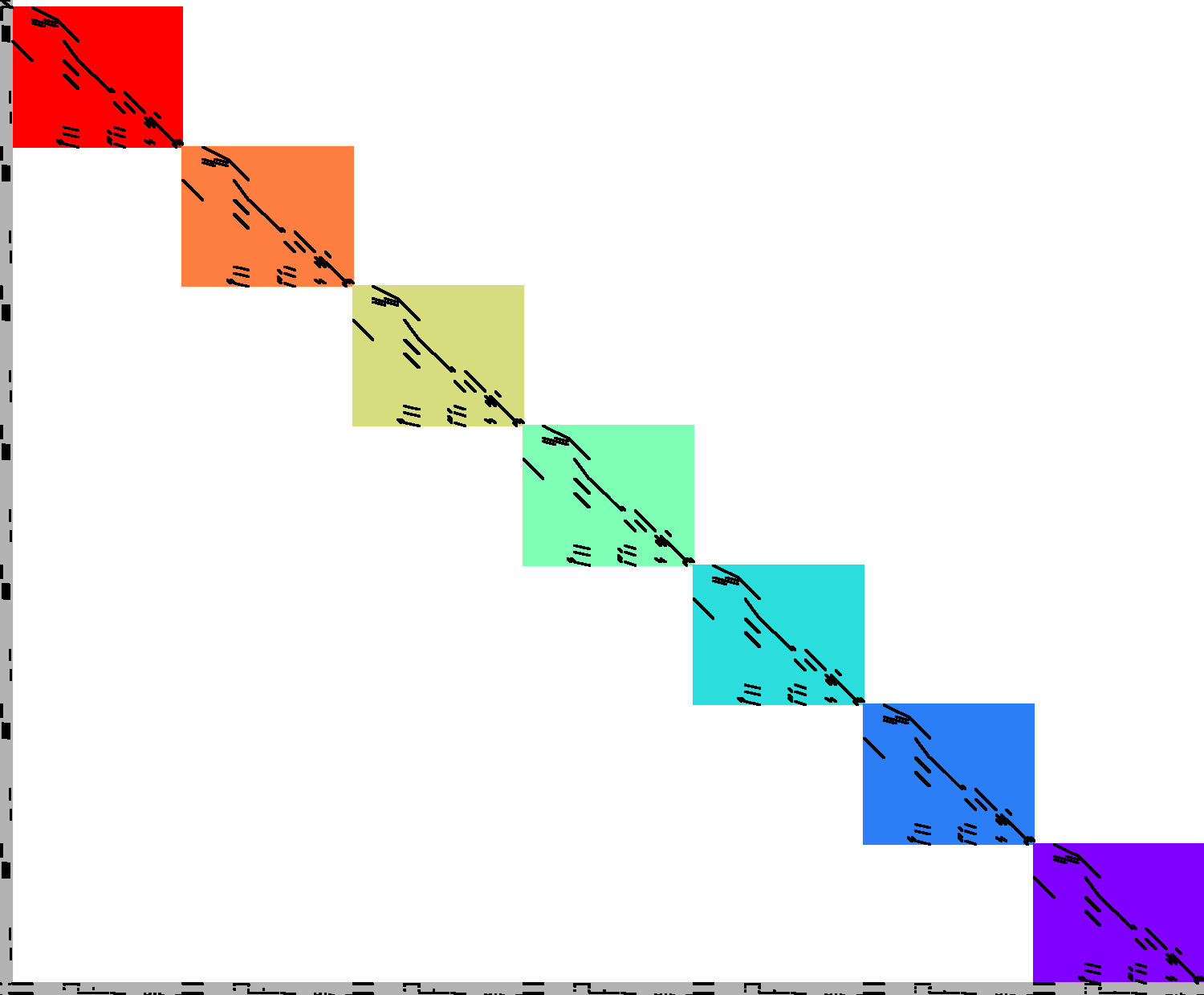}
    \caption{Constraint matrix with the arrowhead structure of a real-world \ESM.}
    \label{fig:remix_blockstructure}
\end{figure}

A recurring structural pattern in linear programs (\LPs)---and more generally, in mixed-integer linear programs (\MIPs)---is the \emph{arrowhead} or \emph{doubly bordered block-diagonal form}, as illustrated in \cref{fig:remix_blockstructure}.
Arrowhead \LPs (\AHLPs) contain two types of linking elements: \emph{linking variables}, which vertically connect multiple diagonal blocks, and \emph{linking constraints}, which connect blocks horizontally.
This structure generalizes both primal and dual block-angular problems.
In many \AHLPs, linking variables and constraints exhibit \emph{local structure}, connecting only a few—often two—consecutive blocks.
This locality arises in a wide range of practical applications.
Energy system models (\ESMs), such as electricity market models with dispatch decisions~\cite{Rehfeldt2019_PIPSIPMpp}, renewable expansion planning and dispatch~\cite{REMIX_2017,PyPSAEur}, and large-scale (stochastic) economic (re-)dispatch models, can often be decomposed spatially or temporally.
Such decompositions naturally give rise to local linking constraints, for example through storage coupling or ramping limits.
Local linking structure also arises in multi-stage \LPs, where decisions at each stage depend on preceding stages.
This includes multi-stage stochastic \LPs, covering applications like asset–liability management, supply network design and supply chain planning, revenue management, and portfolio optimization~\cite{Castro_MultiStage,Colombo2011,Steinbach2001}.
Similarly, staircase \LPs as used in production scheduling, inventory management, transportation, control, and design of multistage structures, exhibit local linking structure connecting consecutive diagonal blocks~\cite{Fourer1982_StaircaseLPs,Wittrock1985_DualNestedOfStaircaseLP}.
Band-diagonal matrices arising in distribution planning can also be reformulated to expose arrowhead forms with local linkage~\cite{OOPS_DistributionPlanning}.
While many optimization problems include integer variables (e.g., for unit commitment or network expansion planning), practitioners often prefer \LP formulations because these scale well computationally and solving \MIPs to optimality remains challenging for large instances, especially in the context of \ESMs~\cite{Anymod,WetzelGreenEnergyCarriers2023,Balmorel}.
Even when a \MIP is ultimately required, its \LP relaxation remains critical, as many heuristic or decomposition-based algorithms rely on (near-)optimal \LP solutions~\cite{MasterThesisTimo_2006}.
This structure is not limited to \LPs and \MIPs.
Examples include nonlinear portfolio optimization~\cite{GondzioOOPS_NonlinPortfolioOptimization}, nonlinear dynamic optimization~\cite{Word2014_EfficientParallelIPMDAT}, model predictive control~\cite{Rao1998_IPMforModelPredictiveControl}, nonlinear parameter estimation~\cite{Zavala2008_IPMForParamEstimation}, and multi-stage nonlinear programs, such as the stochastic optimal power flow model in~\cite{MadNLPGPU}, all of which share the feature of localized linking structure.

Despite decades of progress in general-purpose \LP solvers~\cite{Bixby02, Gondzio12}, they often fail to solve large-scale \AHLPs within reasonable time.
With parallel hardware becoming increasingly available and affordable, one promising direction is the development of highly parallel, structure-exploiting algorithms.
General \LP solution techniques often scale poorly on such hardware, where specialized solutions promise to overcome scaling issues and push back the current computational limits by exploiting the power of high-performance computing (\HPC).
There has been much research activity in decomposing \AHLPs---and nonlinear programs---and their subproblems.
Solution techniques include specialized Simplex methods~\cite{Fourer1982_StaircaseLPs,Friedlander1990_StaircaseLPSimplex}, Danzig-Wolfe decomposition~\cite{Wittrock1985_DualNestedOfStaircaseLP}, Benders-Decomposition~\cite{Meersman2023_NestedBendersStochastic,Zhang2024_BendersForLargeStochasticLP}, Lagrangian decomposition~\cite{kim19}, and frequently interior-point methods (\IPMs)~\cite{CastroIpmForMnf2000,BorderedBlockDiagIPM1996,LubinDualDecomp2013}.
\IPM solvers that specifically exploit this structure include: BlockIp~\cite{Castro2007_BlockIP}, \OOPS~\cite{Gondzio2008_OOPS2}, \PIPSPETRA~\cite{Petra2014_PIPSaugmented}, and \MadNLP~\cite{MadNLPGPU}.

In \IPMs, the focus of this work, often most of the algorithmic time is spent in solving (large-scale) sparse linear saddle-point systems for computing the Newton search direction of the method.
This enables \IPMs to employ specialized linear algebra routines directly working with the structure of a given problem class.

There are two main approaches for solving these systems: sparse direct factorizations (\cite{DuffErismanReid2017_DirectMethodForSparseMatrices}) and iterative schemes with suitable preconditioners (\cite{Saad2003_IterativeMethodForSparseLS}), leading to so-called inexact Newton methods.
We refer the reader to \cite{Gondzio2025_IPMin2025} for an overview of the current state of the art of sparse direct and iterative solvers within \IPMs.
While factorization routines of sparse direct solvers can be slow and memory intensive for large-scale matrices, they often prove robust and exact enough to compute precise search directions.
Additionally, direct factorizations readily support corrector step schemes (\cite{MehrotraPrimalDualIPM1992, gondzio}), solving the same linear system multiple times per iteration at little additional cost.
High precision is especially required towards the end of \IPMs, when the arising linear systems become increasingly ill-conditioned.
During the earlier stages of an \IPM, low accuracy solutions of the underlying linear systems can suffice \cite{Zanetti2023_NewStoppinCriterionForKrylov}.
Iterative schemes, often based on Krylov-subspace methods, are generally more lightweight in their computational effort and memory requirements but struggle to compute high accuracy solutions when the underlying matrix is ill-conditioned.
They heavily rely on strong preconditioners which often need to be recomputed at every (other) iteration of the \IPM.
As saddle-point systems frequently appear in real-world applications their preconditioners have been studied extensively, also in the context of \IPMs (\cite{BenziGolubLiesen_2005_NumericalSaddlePointSystems,Rees2007_PreconditionerForLSFromIPM,DApuzzo2008_MutualImpactOfNumericalLAAndLargeScaleIPM,KarimSolomonik2022_EfficientPreconditionersForIPMViaSCStrategy}).
There also exist notable examples of preconditioners and iterative \IPMs especially exploiting the \AHLP structure.
In \PIPSPETRA and \PIPS (\cite{Petra2014_PIPSaugmented,Rehfeldt2019_PIPSIPMpp}) approximate decomposition methods are combined with Krylov subspace methods to solve large-scale \AHLPs to high accuracy.
In \cite{Lueg2025_DomainDecompositionPreconditioners} the authors suggest several distributed preconditioners based on a distributed Schur complement method underlying an \IPM algorithm. 
Also in the context of pure linear solvers much work has been dedicated to preconditioners exploiting \AHLP matrix structures.
In \cite{SaadSuchomel2002_ARMS} the authors propose the \allcaps{ARMS} linear solver applying a multilevel preconditioner based on incomplete LU factorization to a recursive nested-dissection ordering of their matrices (leading to nested arrowhead structures).
In \cite{XiLiSaad2016_AlgebraicMultilevelPreconditioner} the authors describe a multi-level Schur preconditioner (\allcaps{MLSR}) using a top-down approach approximating the inverse of the original matrix which in \cite{GeoffreyEtAl2018_HierLowRankSCPrecondForIndef} is further refined for indefinite systems, yielding \allcaps{GMLSR}.
A distributed implementation of \allcaps{GMLSR} is discussed in \cite{Xu2022_ParGEMSLR}.
Nevertheless, most commercial optimization software, including \CPLEX~\cite{CPLEX12}, \GUROBI~\cite{Gurobi11}, and the Cardinal Optimizer (\COPT)~\cite{COPT71}, continues to rely on direct (exact) factorizations, valued for their stability and reliability across a wide range of problems.

In this work, we propose a hierarchical direct factorization scheme underlying an \IPM.
In previous work~\cite{Rehfeldt2019_PIPSIPMpp}, we extended the solver \PIPSPETRA, based on \OOQP~\cite{OOQP_WrightGertz}, from primal block-angular problems to \AHLPs.
The resulting solver \PIPS could handle large-scale instances, however, its main bottleneck was the assembly and solution of the so-called (approximate) Schur complement.
The size of the Schur complement is directly tied to the amount of linking structure present in the \AHLP.
Since factorizing and solving with the Schur complement were performed sequentially or in shared-memory parallelism, this posed a bottleneck and limited the scalability of the algorithm.

While iterative or inexact hierarchical schemes based on approximate factorizations can appear attractive for reducing computational effort, our experience indicates that they often suffer from error propagation across recursive elimination levels.
In \IPMs, where the \KKT systems become increasingly ill-conditioned as the iterates approach optimality, such accumulated inexactness can severely impair both numerical stability and the accuracy of the Newton search direction.
Within the original \PIPS framework, we observed that preconditioned BiCGStab could become unstable---even when the preconditioner was recomputed at each \IPM iteration.
For these reasons, we adopt a fully exact hierarchical direct factorization to ensure robustness and reproducibility, while maintaining parallel scalability through distributed elimination.
Nonetheless, we see considerable potential for inexact hierarchical schemes as a complementary direction for future research.

The hierarchical ordering underlying our approach is conceptually related to orderings obtained through independent-set and nested-dissection strategies.
In \cite{Saad1999_BILUM} and subsequently \cite{SaadSuchomel2002_ARMS} the authors use \emph{block independent sets} to split a matrix into arrowhead form and compute an incomplete LU factorization for this system.
They then recursively apply their method to the obtained Schur complement block.
In~\cite{XiLiSaad2016_AlgebraicMultilevelPreconditioner} the authors introduce \allcaps{HID}, the \emph{Hierarchical Interface Decomposition}, based on nested dissection (\allcaps{ND}) which generates matrices that can be reordered into hierarchical arrowhead form (though their focus also lies on the Schur complement structure).
Popular graph partitioning packages such as \allcaps{hMETIS} \cite{Karypis1997_HMETIS}, PaToH \cite{Catalyuerek2011_PaToH}, and KaHyPar \cite{Schlag2023_KaHyPar} based on multi-level hyper-graph partitioning can also be used to obtain similar matrix reorderings.
The generic decomposition solver \allcaps{GCG}~\cite{Gamrath2010_GenericDantzigWolfeGCG} includes an arrowhead detection method based on \allcaps{hMETIS}.

While these methods might be used to obtain an ordering generally suitable for \HSCA, the main challenge lies in generating a well balanced block structure at each level of the algorithm.
In \HSCA, the arrowhead structure detection happens already on the \LP system matrix rather than the \KKT matrix.
This allows us to employ a complex suite of structure preserving presolving techniques \cite{Kempke19} (an extremely important part of solving an \LP) to the \LP formulation of our problem, without destroying the arrowhead structure or its well balanced-ness.
Currently, \HSCA relies on modeler-supplied decompositions that reflect the underlying physical or temporal structure, ensuring extremely well-balanced diagonal blocks and explicit control over the hierarchy of local links.
However, we view automatic structure detection as an important part of a structure exploiting solver and are experimenting with detection mechanisms for large-scale matrices based on the KarHyPar package\footnote{https://gitlab.com/pips-ipmpp/detection-annotation}.

The main contribution of this article is the description and implementation of a novel approach within \PIPS, called the \emph{hierarchical Schur complement approach} (\HSCA).
Exploiting the local structure of linking constraints, we recursively decompose the Schur complement into smaller, distributed components that can be solved in parallel.
While our implementation specifically targets linking constraints connecting two consecutive blocks, most prevalent in our models, constraints connecting multiple adjacent blocks can be treated analogously and the locality of linking variables can be similarly exploited.
As \IPMs also form the foundation for many quadratic and nonlinear solvers, the proposed method extends naturally to a wide class of structured nonlinear programs.

Our implementation is motivated by large-scale economic dispatch models with storage expansion.
For the computational evaluation of \HSCA, we rely on two different types of \ESMs.
First, we tested on the \SIMPLE instances developed during the research projects \beamme\footnote{\beamme project: \url{http://www.beam-me-projekt.de/beam-me/EN/Home/home_node.html}} and \unseen\footnote{\UNSEEN project: \url{https://unseen-project.gitlab.io/home/}}, see~\cite{Rehfeldt2019_PIPSIPMpp,best_practice_guide}.
These simplified models capture the essential features of \ESMs, including all necessary DC power flow and storage balance constraints.
The number and size of the diagonal blocks can be precisely controlled without incurring additional modeling overhead.
Second, we consider a family of economic redispatch models generated with \remixmiso, a variant of the \remix framework~\cite{REMIX_2017} developed at the German Aerospace Center\footnote{German Aerospace Center: \url{https://www.dlr.de/de}}.
These models are described in detail in~\cite{yssp,Wetzel2024_REMIX}.
Using \HSCA, we successfully solve instances with up to 1.7 million linking constraints, demonstrating substantial improvements in both scalability and performance.

\section{Interior-point methods for linear programming}

The two main methods for solving general \LPs are the simplex algorithm~\cite{chvatal1983_LP} and \IPMs~\cite{wright1997_primal_dual_ipm}.
Often, \IPMs are more successful on large problems and offer more potential for parallelization.
This success can be explained by the fact that the \IPM's main computational effort is matrix factorization, and the amount of factorizations required grows relatively slowly with the problem size~\cite{PotraWrightIPMs}.
Many algorithms fall into the notion of \IPMs, see \cite{terlaky2013_IPM, wright1997_primal_dual_ipm}.
In this work, we use an infeasible primal-dual \IPM similar to the ones described in \cite{LustigMehrotra1992,MehrotraPrimalDualIPM1992} using Mehrotra's predictor corrector scheme extended by Gondzio corrector steps \cite{gondzio}.
In the following, we omit many details on the \IPM method, and instead focus on the main feature of \HSCA: the hierarchical decomposition of \AHLPs.
Details of the actual \IPM can be found in \cite{OOQP_WrightGertz,LubinEtAl2011,Petra2014_PIPSaugmented} and the references mentioned therein.

\subsection{Notation}

We often write $0$ for a zero-vector or matrix.
As is often the case, we omit zero entries for large matrices with more complicated block structures.
We sometimes selectively write $0$ and omit entries, whenever it helps depict a given matrix block structure.
$I_n$ denotes the identity matrix of order $n$.
When it is clear from the context, we omit the index and write $I$.
The vector consisting of all ones, e.g., $(1,\dots,1)\T\in\R^n$, will be denoted by $\vecOnes$.
For a vector $x \in \R^n$, $X = \diag(x) \in \R^{n \times n}$ denotes the diagonal matrix with $x$ on its diagonal.
For a given matrix $A \in \R^{m \times n}$ and two subsets $I \subset \{1, \dots, m\}$, $J \subset\{1, \dots, n\}$, we refer to the submatrix built from $A$ by taking only the rows/columns indicated by $I/J$ as $A_{I,.}$ /$A_{.,J}$ respectively.

\subsection{Linear programming, duality and optimality conditions}

Consider the following \LP in standard form
\begin{alignat*}{2}
    \min_x \quad& \objectiveLinearVector\T x \\
    \subjectto \quad&\coefmatrixEq x = b \\
    & x \ge 0,
\end{alignat*}
where $c \in \R^n, b\in\R^m, A \in \Rmn$, as well as its associated dual formulation
\begin{alignat*}{2}
    \max_{y, \dualsLowerBounds} \quad& b\T y \\
    \subjectto \quad&\coefmatrixEq \T y + \dualsLowerBounds = \objectiveLinearVector \\
    & \dualsLowerBounds \ge 0,
\end{alignat*}
where $y \in \R^m, \dualsLowerBounds \in \R^n$.
We assume that all matrices considered here have full rank.
It is well-known that if both problems are feasible, the optimal primal and dual objectives coincide \cite{wright1997_primal_dual_ipm}.
The necessary \KKT conditions to either the primal or the dual problem (which motivates the term primal-dual) are given by:
\begin{align}
\coefmatrixEq\T \dualsEq + \dualsLowerBounds - \objectiveLinearVector& = 0 \label{eq:kkt_stationarity} \\
\coefmatrixEq \primalvector - \rhsvectorEq & = 0 \label{eq:kkt_primal_feasibility} \\
\dualsLowerBoundsCap \primalvectorCap \vecOnes &= 0 \label{eq:kkt_complementarity} \\
\primalvector, \dualsLowerBounds & \ge 0, \label{eq:kkt_feasibility}
\end{align}
where $X = \diag(x)$ and $\dualsLowerBoundsCap = \diag(\dualsLowerBounds)$.
All but the complementarity equation \cref{eq:kkt_complementarity} are linear.

\subsection{Introduction to primal-dual interior-point methods}

In primal-dual infeasible \IPMs, perturbed \KKT conditions are solved by Newton's method, namely \cref{eq:kkt_complementarity} is perturbed by $\centralityFactor > 0$:
\begin{equation}
\dualsLowerBoundsCap \primalvectorCap \vecOnes = \centralityFactor \vecOnes.
\label{eq:kkt_perturbed_complementarity}
\end{equation}
The so-called centrality factor $\centralityFactor$ is successively decreased, which guides the \IPM along the \textit{central path} defined by the values of $\centralityFactor$ close to the optimal solution until sufficient convergence has been reached.
In doing so, the \IPM maintains its iterates $\primalvector, \dualsLowerBounds$ strictly within the positive orthant.

We introduce the vector notation for \cref{eq:kkt_stationarity,eq:kkt_primal_feasibility,eq:kkt_feasibility,eq:kkt_perturbed_complementarity} as well as an iteration index $k$.
At each step $k$ of the \IPM, we approximately solve
\begin{equation}\label{eq:newton_system}
\matr{
\coefmatrixEq\T \dualsEq^{k} + \dualsLowerBounds^{k} \\ \coefmatrixEq \primalvector^{k} \\
\primalvectorCap^{k} \dualsLowerBoundsCap^{k} \vecOnes
}
=
\matr{
\objectiveLinearVector \\
\rhsvectorEq \\
\centralityFactor^{k} \vecOnes
}.
\end{equation}
%
%
%
The linear system that arises in each iteration $l$ of Newton's method is
\begin{equation*}
\matr{
    0 & \coefmatrixEq\T & I \\
    \coefmatrixEq & 0 & 0   \\
    \dualsLowerBoundsCap^{k,l} & 0 & \primalvectorCap^{k,l}
    }
    \matr{\D\primalvector^{k,l}\\ \D\dualsEq^{k,l} \\ \D\dualsLowerBounds^{k,l}} = \matr{ r_{\primalvector^{k,l}} \\ r_{\dualsEq^{k,l}} \\ r_{\dualsLowerBounds^{k,l}}},
\end{equation*}
where $r_{\primalvector^{k,l}}, r_{\dualsEq^{k,l}}, r_{\dualsLowerBounds^{k,l}}$ are residual right-hand sides.
Recalling that $\primalvector^{k,l} > 0$, we eliminate the last row block
\begin{equation*}
    \D\dualsLowerBounds^{k,l} = (\primalvectorCap^{k,l})^{-1} r_{\dualsLowerBounds^{k,l}} - (\primalvectorCap^{k,l})^{-1}\dualsLowerBoundsCap^{k,l} \D\primalvector^{k,l},
\end{equation*}
which produces the symmetric so-called \textit{augmented system}~\cite{MarosMeszarosAugmentedSystem1998, wrightAugSys1997}
\begin{equation}\label{eq:linear_system_aug_ipm}
\matr{
    -(\primalvectorCap^{k,l})^{-1} \dualsLowerBoundsCap^{k,l} & \coefmatrixEq\T \\
    \coefmatrixEq & 0
} \matr{\D\primalvector^{k,l}\\ \D\dualsEq^{k,l}}= \matr{ \hat{r}_{\primalvector^{k,l}} \\ r_{\dualsEq^{k,l}} }.
\end{equation}
In \IPMs for \LP, the number of inner iterations for approximately solving \cref{eq:newton_system} is generally kept to a single inner iteration $l$, after which the centrality factor $\tau$ is updated.

\subsection{The Schur complement algorithm}

Here, we recall the Schur complement algorithm, slightly modified for our purposes.
We will then use it to derive \HSCA.
It can be thought of as a block Gaussian elimination.
Consider the linear system
\begin{equation}\label{eq:generalized_sc_system}
    \matr{ K & L \\ L^T & K_0 } \matr{ x_1 \\ x_0 } = \matr{ b_1 \\ b_0 }
\end{equation}
with matrices $K \in \R^{m\times m}$, $L \in\R^{m\times n}$, $K_0 \in \R^{n \times n}$, $n$, $m\in\N$ and vectors $x_1$, $x_0$, $b_1$, $b_0$ of appropriate dimensions.
Given a factorization of $K$, we can compute and factorize the Schur complement $S := K_0 - L\T K\inv L$ of \cref{eq:generalized_sc_system} column-wise, as described in \cref{alg:schur_complement_general}, and solve the linear system accordingly.
The factorization of $K$ is not required to be given as actual factors of a matrix factorization.
Instead, any algorithm that allows solving linear equations with $K$ suffices for the computation of $S$.
We refer to such an algorithm as an \textit{implicit factorization} of $K$.
\begin{algorithm}
    \caption{Schur Complement Factorization}
    \label{alg:schur_complement_general}
    \hspace*{\algorithmicindent}~~\textbf{Input:} System \cref{eq:generalized_sc_system}
    \begin{algorithmic}[1]
        \setcounter{ALC@unique}{0}
        \STATE{Set $S=K_0$}
        \FOR{$i\in\{1,\dots,n\}$}
            \STATE Solve $K z = L_{.,i}$
            \STATE Set $S_{.,i} = S_{.,i} - L^T z$
        \ENDFOR
        \STATE{Factorize $S$}
        \RETURN
    \end{algorithmic}
\end{algorithm}
\section{A specialized parallel interior-point method}
\label{sec:pips_ipm}

In the remainder of this section, we will show how \PIPS employs a parallel Schur complement decomposition~\cite{Rehfeldt2019_PIPSIPMpp} to solve \cref{eq:linear_system_aug_ipm} in a distributed fashion and in parallel.
This forms the basis for the methods presented in \cref{sec:hierarchical_approach}.

\subsection{Message passing interface}

We will be relying on the availability of communication between independent processes on a distributed computing system.
Our implementation heavily uses the de facto standard in messaging protocols for distributed-memory computing, the message passing interface \MPI \cite{LydonEtAl1994_MPIStandard}.
\MPI provides a set of collective and point-to-point communication routines for a given set of parallel processes.
We commonly refer to these processes connected via \MPI as \textit{\MPI processes} or just \textit{processes}.
Each of the \MPI processes is assigned a unique global rank in ${0, \dots, N-1}$, referred to as \textit{global rank}, where $N$ is the total number of processes.
The choice of $N$ as both the number of blocks and the number of \MPI processes is intended, as they are essentially the same (as explained in the following chapters).
The global ranks are associated with a global \MPI communicator.
Generally, communicators are used by processes for collective communication operations such as \textit{Reduce} and \textit{Allreduce}, where a specific reduction operation, e.g., summation, is used to combine the data of all \MPI processes on respectively one or all processes.
\MPI offers the option of creating custom sub-communicators, each of which can contain an arbitrary subset of the $N$ processes and assigns an additional (communicator-unique) rank to its processes.
When discussing a process's rank, we mean its associated global rank (exceptions will be stated explicitly).
Throughout this work, we use the Allreduce and Reduce operations to operate on vectors and matrices.
Additionally, we use the concept of sub-communicators for the distribution of the linear algebra in \HSCA.

\subsection{A parallel Schur complement decomposition}
\label{sec:parallel_sc_decomp}

The \AHLPs whose system matrix is depicted in \cref{fig:remix_blockstructure} can mathematically be described as

\begin{equation}\label{fig:blockstructure}
\begin{alignedat}{4}
    \min \quad {{c_0^T x_0} }~~+~~			& { c_1^T x_1 }~~+~~\cdots		&~~+~~c_N^T x_N		& \\
    \subjectto \quad {\coefmatrixEq_0 x_0} ~~~~~~~~ 			&   								 	 		&										& = {b_0} \\
    {d_0}	\leq{\coefmatrixIneq_0  {x_0}} ~~~~~~~~ 	&  									  			&										&\leq{f_0}  \\
    {\coefmatrixEq_1 {x_0} }~~+~~     						&  B_1 x_1 		 		        &        								&= b_1  \\
    { d_1}\leq{  \coefmatrixIneq_1  {x_0}}~~+~~& D_1 x_1 			            &        								&\leq f_1  \\
    { \vdots}~~~~~~~~~~											&~~~~~~~~~~~~~~\ddots     &   	  														&~\vdots  \\
    { 	\coefmatrixEq_N{x_0}} ~~+~~      					&        							            &~~+~~ B_N x_N	  		&= b_N  \\
    { 	d_N	}\leq{\coefmatrixIneq_N {x_0}}~~+~~&				                    	       	&~~+~~ D_N x_N 		&\leq f_N  \\
    { {	{\PIPSBl_0 x_0}} }~~+~~  &{  \PIPSBl_1 x_1}~~+~~\cdots &~~+~~\PIPSBl_N x_N      	&= {b_{N+1}}  \\
    { {d_{N+1}}}\leq{{\PIPSDl_0 x_0}}~~+~~   &{  \PIPSDl_1 x_1}~~+~~\cdots &~~+~~\PIPSDl_N x_N			& \leq {f_{N+1}} \\
    && \ell_i \leq x_i &\leq u_i \quad \forall i=0,\dots,N.
\end{alignedat}
\end{equation}

Here $x_i \in \R^{n_i}$ denote the decision variables, $c_i \in \R^{n_i}$ the objective vector and $\ell_i$, $u_i\in\R^{n_i}$ the variable bounds.
The vectors $b_i\in\R^{m_{i_\coefmatrixEq}}$, $d_i$, and $f_i\in\R^{m_{i_\coefmatrixIneq}}$ denote the right-hand sides for equalities and lower and upper bounds for inequalities respectively.
The system matrix is split into the sub-matrices $A_i \in \R^{m_{i_\coefmatrixEq}\times n_0}$, $B_i\in\R^{m_{i_\coefmatrixEq}\times n_i}$, $\coefmatrixIneq_i\in\R^{m_{i_\coefmatrixIneq}\times n_0}$, $D_i\in\R^{m_{i_\coefmatrixIneq}\times n_0}$, $\PIPSBl_i \in \R^{m_{{N+1}_\coefmatrixEq}\times n_i}$, $\PIPSDl_i \in \R^{m_{{N+1}_\coefmatrixIneq}\times n_i}$ for $N \in \N$.
We denote by $n_i\in \N$, the size of the $i$-th primal vector, by $m_{i_\coefmatrixEq} \in \N$ the row dimension of the $i$-th equality constraint sub-matrix $A_i$ and by $m_{i_\coefmatrixIneq}\in\N_0$ the row dimension of the $i$-th inequality constraint sub-matrix $C_i$.
Using this notation, linking constraints are constraints associated with the blocks $\PIPSBl_i$ and $\PIPSDl_i$, while linking variables are the variables in $x_0$.
In the following, we will further distinguish linking constraints into \textit{local linking constraints} and \textit{global linking constraints}.
Local linking constraints (also called two-links) only contain nonzeros in $\PIPSBl_0$ (or $\PIPSDl_0$) and two consecutive blocks $\PIPSBl_i$ and $\PIPSBl_{i+1}$ (or $\PIPSDl_i$ and $\PIPSDl_{i+1})$ for exactly one $i \in \{1,\dots,N-1\}$.
Global linking constraints constitute the rest of the linking constraints, in particular the rows that have entries in $\PIPSBl_0, \PIPSBl_i, \PIPSBl_j$ (or $\PIPSDl_0, \PIPSDl_i, \PIPSDl_j$) where $|i-j| > 1$ and $i,j \in \{1,\dots,N\}$.

For clarity, we will omit the outer and inner iteration indices $k, l$ introduced for \IPMs in the rest of this section.
Similarly, we will omit the $\coefmatrixIneq$, $D$, and $G$ matrices associated with inequality constraints of system \cref{fig:blockstructure}.
Consequently, we refer to the sizes of the remaining $A_i$ as $\R^{m_i\times n_0}$, $B_i$ as $\R^{m_i\times n_i}$ and $F_i$ as $\R^{m_{N+1}\times n_i}$, dropping the matrix index in $m_{i_A}$.

The augmented system \cref{eq:linear_system_aug_ipm} corresponding to \cref{fig:blockstructure} is
\begin{equation*}
    \matr{
        \Sigma_0 & & & & A_0\T & A_1\T & \dots  & A_N\T & \PIPSBl_0\T \\
        & \Sigma_1 & & &       & B_1\T &        &       & \PIPSBl_1\T \\
        & &\ddots&&       &       & \ddots &       & \vdots      \\
        & & & \Sigma_N &       &       &        & B_N\T & \PIPSBl_N\T \\
        A_0 & & & & & & & & \\
        A_1 & B_1 & & & & & & & \\
        \vdots & & \ddots & & & & & & \\
        A_N & & & B_N & & & & & \\
        \PIPSBl_{0} & \PIPSBl_1 & \dots & \PIPSBl_N & & & & &
    }
    \matr{
        \Delta \primalvector_0 \\
        \Delta \primalvector_1 \\
        \vdots \\
        \Delta \primalvector_N \\
        \Delta \dualsEq_0 \\
        \Delta \dualsEq_1 \\
        \vdots \\
        \Delta \dualsEq_N \\
        \Delta \dualsEq_{N+1}
        } =
    \matr{ r_{\primalvector_0} \\ r_{\primalvector_1} \\ \vdots \\ r_{\primalvector_N} \\ r_{\dualsEq_0} \\ r_{\dualsEq_1} \\ \vdots \\ r_{\dualsEq_N} \\ r_{\dualsEq_{N+1} } },
\end{equation*}
where
$\Sigma_i := -\primalvectorCap_i\inv\dualsLowerBoundsCap_i$ for $i={0, \dots, N}$.
After a symmetric permutation, we get
%
%
\begin{equation}\label{eq:linear_system_pips_perturbed_compressed}
\matr{
    K_1 & & & L_1 \\
    & \ddots & & \vdots \\
    & & K_N & L_N \\
    L_1\T & \dots & L_N\T & K_0
} \matr{ \D z_1 \\ \vdots \\ \D z_N \\ \D z_0 } = \matr{ b_1 \\ \vdots \\ b_N \\ b_0 },
\end{equation}
where
\begin{equation*}
K_i := \matr{
\Sigma_i & B_i\T \\ B_i & 0
},\
K_0 := \matr{
\Sigma_0 & A_0\T & \PIPSBl_0\T \\
A_0 & 0 & 0 \\
\PIPSBl_0 & 0 & 0
},\
L_i := \matr{
0 & 0 & \PIPSBl_i\T \\
A_i & 0 & 0
},
\end{equation*}
\begin{equation*}
\D z_0 := \matr{
    \D \primalvector_0 \\
    \D \dualsEq_0 \\
    \D \dualsEq_{N+1}
},\
\D z_i := \matr{
    \D \primalvector_i \\
    \D \dualsEq_i
},\
b_0 := \matr{
    r_{\primalvector_0} \\
    r_{\dualsEq_0} \\
    r_{\dualsEq_{N+1}}
},\
b_i := \matr{
    r_{\primalvector_i} \\
    r_{\dualsEq_i} \\
},
\end{equation*}
for $i=\{1,\dots,N\}$.
As described in \cite{Rehfeldt2019_PIPSIPMpp}, this system can be solved efficiently in parallel using a Schur complement decomposition.

As usual in direct methods, we distinguish between two phases, \textit{factor} and \textit{solve}.
The factorization phase for solving system \cref{eq:linear_system_pips_perturbed_compressed} is described by \cref{alg:factor_pips}.
The solve phase for a given right-hand side is described by \cref{alg:solve_pips}.

\begin{algorithm}
    \caption{{\textit{Parallel}} factorization}
    \label{alg:factor_pips}
    \hspace*{\algorithmicindent}~~\textbf{Input:} System matrix \cref{eq:linear_system_pips_perturbed_compressed}
    \begin{algorithmic}[1]
        \setcounter{ALC@unique}{0}
        \STATE\label{alg:factor_pips:factorizeK}{\textit{Factorize (potentially implicit)} $K_i$ for all $i\in\{1,\dots,N\}$}
        \STATE\label{alg:factor_pips:computeSCcontrib}{\textit{Compute} $-L_i\T K_i\inv L_i$ for all $i\in\{1,\dots,N\}$}
        \STATE\label{alg:factor_pips:SumSC}{Sum Schur complement $S := K_0 - \sum_{i=1}^{N} L_i\T K_i\inv L_i$}
        \STATE\label{alg:factor_pips:FactorSC}{Factorize $S$}
        \RETURN
    \end{algorithmic}
\end{algorithm}

\begin{algorithm}
    \caption{{\textit{Parallel}} solve 
    }
    \label{alg:solve_pips}
    \hspace*{\algorithmicindent}~~\textbf{Input:} Factorized system from \cref{alg:factor_pips} and right hand side in \cref{eq:linear_system_pips_perturbed_compressed}
    \\ \hspace*{\algorithmicindent}~~\textbf{Output: } Solution $\D z_i,\ i\in\{0,\dots,N\}$
    \begin{algorithmic}[1]
        \setcounter{ALC@unique}{0}
        \STATE\label{alg:solve_pips:computeSCrhsContrib}{\textit{Compute} $-L_i\T K_i\inv b_i$ for all $i\in\{1,\dots,N\}$}
        \STATE\label{alg:solve_pips:SumSCrhs}{Sum Schur complement right hand side $\hat{b}_0 = b_0 - \sum_{i=1}^{N} L_i\T K_i\inv b_i$}
        \STATE\label{alg:solve_pips:solveSCrhs}{Solve $S \D z_0 = \hat{b}_0$}
        \STATE\label{alg:solve_pips:computeModRhs}{\textit{Compute modified right-hand sides} $\hat{b}_i = b_i - L_i \D z_0$}
        \STATE\label{alg:solve_pips:solveZi}{\textit{Solve} $K_i \D z_i = \hat{b}_i$ for all $i\in\{1,\dots,N\}$}
        \RETURN
    \end{algorithmic}
\end{algorithm}
\Cref{alg:factor_pips} extends \cref{alg:schur_complement_general} in \PIPS.
In its first implementation, \PIPS used a combination of \cref{alg:schur_complement_general} and \cref{alg:factor_pips} to compute the Schur complement.
Instead of forming $K_i\inv$ explicitly, the Schur complement contributions $L_i\T K_i\inv L_i$ were computed by factoring $K_i$ and solving $K_i z = l$ for each column $l$ of $L_i$.
Currently, the faster way of computing the blockwise Schur complement contributions $L_i K_i^{-1} L_i\T$ is the application of an incomplete (partial) LU factorization routine~\cite{Petra2014_PIPSaugmented} to the extended system
\begin{equation*}
\matr{K_i & L_i\T \\ L_i & 0}.
\end{equation*}
By restricting pivoting to the $(1,1)$ block and aborting after $m_i + n_i$ pivots, the transformed $(2,2)$ block contains $-L_i\T K_i^{-1} L_i$.
This approach is implemented in the solver PARDISO~\cite{pardiso-7.2a}.
Similar to \Cref{alg:schur_complement_general} is it not necessary to compute actual factors for the factorization of $K_i$ but it suffices to have an implicit factorization available.

\cref{alg:factor_pips:factorizeK,alg:factor_pips:computeSCcontrib} of \cref{alg:factor_pips}, and \cref{alg:solve_pips:computeSCrhsContrib,alg:solve_pips:computeModRhs,alg:solve_pips:solveZi} of \cref{alg:solve_pips} can be executed in parallel and are highlighted in italic.
When each \MPI process has access to the data needed for the parallel steps, everything but the Schur complement factorization and the Schur complement solve can be run in distributed parallel.
No process needs a representation of the whole problem: a process requires knowledge of only the system blocks $i$ assigned to it and the $0$ block as depicted in \cref{fig:blockstructure}.
This also induces a limit on the number of processes and motivates the slightly ambiguous notation of $N$ for the number of blocks and processes.
A problem with $N$ blocks can be distributed among at most $N$ processes.
We point out that \cref{alg:factor_pips:SumSC} in \cref{alg:factor_pips} and \cref{alg:solve_pips:SumSCrhs} in \cref{alg:solve_pips} require communication between the processes to form the Schur complement and its right-hand side.
To save memory and compute power, the Schur complement in \PIPS is stored on a single process, rendering the communication in both steps a \textit{reduce operation}.
%
%
Processing the Schur complement on a single \MPI process requires additional communication (\textit{broadcast}) of the solution $z_0$ in \cref{alg:solve_pips:solveSCrhs} of \cref{alg:solve_pips}.
%

\subsection{Detecting structure in the linking constraints}\label{sec:linkcons_structure}

One of the main contributions of \cite{Rehfeldt2019_PIPSIPMpp} was detecting and exploiting additional structure in the linking part of the system matrix.
This will also be essential in \cref{sec:hierarchical_approach}.
In the following, we summarize the most important findings.

We initially noticed that many linking constraints in \ESMs were indeed local linking constraints (they link exactly two consecutive blocks in addition to the $0$th block).
This introduces additional structure to the problem:
we can split $\PIPSBl_i$ into local $\PIPSBlloc_i$ and global $\PIPSBlglob_i$ linking constraints, by defining
\begin{equation*}
\PIPSBl_1 := \matr{
    \PIPSBlloc_1       \\
    0                  \\
    \\
    \vdots             \\
    \\
    \\
    0                  \\
    \PIPSBlglob_1      \\
}, \PIPSBl_2 := \matr{
    \PIPSBlloc'_1      \\
    \PIPSBlloc_2       \\
    0                  \\
    \\
    \vdots             \\
    \\
    0                  \\
    \PIPSBlglob_2      \\
}, \PIPSBl_3 := \matr{
    0                  \\
    \PIPSBlloc'_2      \\
    \PIPSBlloc_3       \\
    0                  \\
    \vdots             \\
    \\
    0                  \\
    \PIPSBlglob_3      \\
}, \dots,\PIPSBl_{N-1} := \matr{
    0                  \\
    \\
    \vdots             \\
    \\
    0                  \\
    \PIPSBlloc_{N-2}'  \\
    \PIPSBlloc_{N-1}   \\
    \PIPSBlglob_{N-1}  \\
}, \PIPSBl_{N} := \matr{
    0                  \\
    \\
    \vdots             \\
    \\
    \\
    0                  \\
    \PIPSBlloc'_{N-1}  \\
    \PIPSBlglob_{N}    \\
},
\end{equation*}
with $\PIPSBlloc_i \in \R^{l_i\times n_i}$, $\PIPSBlloc'_i \in \R^{l_i\times n_{i+1}}$, $l_i \in \N_0$ for $i\in \{1,\dots,N-1\}$ and $\PIPSBlglob_i \in \R^{m_\PIPSBlglob\times n_i}$ with $m_\PIPSBlglob \in \N_0$.
We set
\begin{equation*}
\hat{\PIPSBl}_i := \matr{
0\\
\vdots \\
0\\
\PIPSBlloc'_{i-1} \\
\PIPSBlloc_i \\
0\\
\vdots \\
0\\
},\quad
{\PIPSBl_{0}} := \matr{
    \PIPSBlloc_{0,1}   \\
    \PIPSBlloc_{0,2}   \\
    \\
    \vdots             \\
    \\
    \PIPSBlloc_{0,N-2} \\
    \PIPSBlloc_{0,N-1} \\
    \PIPSBlglob_{0}   \\
}, \quad \hat{\PIPSBl_0} := \matr{
    \PIPSBlloc_{0,1}   \\
    \PIPSBlloc_{0,2}   \\
    \\
    \vdots             \\
    \\
    \PIPSBlloc_{0,N-2} \\
    \PIPSBlloc_{0,N-1} \\
},
\end{equation*}
with $\PIPSBlglob_{0}\in \R^{m_{\PIPSBlglob_0}\times n_0}$, $\PIPSBlloc_{0,i} \in \R^{l_i\times n_0}$, for $i\in\{1,\dots,N-1\}$. Analogously expanding right hand side and variable vectors, the linear system \cref{eq:linear_system_pips_perturbed_compressed} becomes
\begin{equation}
\label{eq:linear_system_pips_perturbed_full_local_global_links}
\left[
\begin{array}{ccccc|cccc}
\Sigma_1        & B_1\T &        &                 &       & 0               & 0        & \hat{\PIPSBl}_1\T    & \PIPSBlglob_1\T    \\
B_1             & 0     &        &                 &       & A_1             & 0        & 0                    & 0                  \\
                &       & \ddots &                 &       &                 &  \vdots  &                &                    \\
                &       &        & \Sigma_N        & B_N\T & 0               & 0        & \hat{\PIPSBl}_N\T    & \PIPSBlglob_N\T    \\
                &       &        & B_N             & 0     & A_N             & 0        & 0                    & 0                  \\
\hline
0               & A_1\T &        & 0               & A_N^T & \Sigma_0        & A_0\T    & \hat{\PIPSBl}_0\T    & \PIPSBlglob_0\T    \\
0               & 0     & \dots  & 0               & 0     & A_0             & 0        & 0                    & 0                  \\
\hat{\PIPSBl}_1 & 0     &        & \hat{\PIPSBl}_N & 0     & \hat{\PIPSBl}_0 & 0        & 0                    & 0                  \\
\PIPSBlglob_1   & 0     &        & \PIPSBlglob_N   & 0     & \PIPSBlglob_0   & 0        & 0                    & 0                  \\
\end{array}
\right] \matr{ \D\primalvector_1 \\ \D\dualsEq_1 \\ \vdots \\ \D\primalvector_N \\ \D\dualsEq_N \\ \D\primalvector_0 \\ \D\dualsEq_0 \\ \D\dualsEq_{N+1} \\\D\mathbf{\dualsEq_{N+1}} } =
\matr{ r_{\primalvector_1} \\ r_{\dualsEq_1} \\ \vdots \\ r_{\primalvector_N} \\ r_{\dualsEq_N} \\ r_{\primalvector_0} \\ r_{\dualsEq_0} \\ r_{\dualsEq_{N+1}} \\ \mathbf{r_{\dualsEq_{N+1}}} }.
\end{equation}
We note that the locality of many of the linking constraints also induces a arrowhead block structure on the Schur complement of \cref{eq:linear_system_pips_perturbed_compressed}.
\PIPS exploits this structure to use a sparse direct solver for the Schur complement.
We refer the reader to \cref{apdx:sc_structure} for more details.
\section{A hierarchical solution approach}
\label{sec:hierarchical_approach}

As already mentioned in \cref{sec:intro}, one of the main bottlenecks in \PIPS for solving large problems or problems with many linking constraints is the size of the Schur complement.
In the following, we outline an algorithm that circumvents this bottleneck by further exploiting the local linking constraints.
It splits the Schur complement into multiple smaller Schur complements hierarchically connected in the factorization and solve processes.

\subsection{The dense layer}

After another symmetric permutation in the borders, which moves the local linking constraints $\hat{\PIPSBl}_i\T$/$\hat{\PIPSBl}_i$ $i=0,\dots,N$ to the front/top of their respective blocks,
\cref{eq:linear_system_pips_perturbed_full_local_global_links} becomes
\begin{equation*}
\left[
\begin{array}{ccccc|c|ccc}
    \Sigma_1        & B_1\T &        &                 &       & \hat{\PIPSBl}_1\T    & 0               & 0     & \PIPSBlglob_1\T   \\
    B_1             & 0     &        &                 &       & 0                    & A_1             & 0     & 0                 \\
                    &       & \ddots &                 &       & \vdots               &                 & \vdots&                   \\
                    &       &        & \Sigma_N        & B_N\T & \hat{\PIPSBl}_N\T    & 0               & 0     & \PIPSBlglob_N\T   \\
                    &       &        & B_N             & 0     & 0                    & A_N             & 0     & 0                 \\
    \hline
    \hat{\PIPSBl}_1 & 0     & \dots  & \hat{\PIPSBl}_N & 0     & 0                    & \hat{\PIPSBl}_0 & 0     & 0                 \\
    \hline
    0               & A_1\T &        & 0               & A_N^T & \hat{\PIPSBl}_0\T    & \Sigma_0        & A_0\T & \PIPSBlglob_0\T   \\
    0               & 0     & \dots  & 0               & 0     & 0                    & A_0             & 0     & 0                 \\
    \PIPSBlglob_1   & 0     &        & \PIPSBlglob_N   & 0     & 0                    & \PIPSBlglob_0   & 0     & 0                 \\
\end{array}
\right] \matr{ \D\primalvector_1 \\ \D\dualsEq_1 \\ \vdots \\ \D\primalvector_N \\ \D\dualsEq_N  \\ \D\dualsEq_{N+1} \\ \D\primalvector_0\\ \D\dualsEq_0\\ \D\mathbf{\dualsEq_{N+1}} } =
\matr{ r_{\primalvector_1} \\ r_{\dualsEq_1} \\ \vdots \\ r_{\primalvector_N} \\ r_{\dualsEq_N} \\ r_{\dualsEq_{N+1}}  \\ r_{\primalvector_0} \\ r_{\dualsEq_0} \\  \mathbf{r_{\dualsEq_{N+1}}}}.
\end{equation*}

We can apply the Schur complement decomposition in \cref{alg:schur_complement_general} by setting
\begin{equation*}
K^{\text{dense}} := \scalemath{0.9}{\matr{
	\Sigma_1     & B_1\T &        &                 &       & \hat{\PIPSBl}_1\T \\
	B_1          & 0     &        &                 &       & 0                 \\
	             &       & \ddots &                 &       & \vdots            \\
	             &       &        & \Sigma_N        & B_N\T & \hat{\PIPSBl}_N\T \\
                 &       &        & B_N             & 0     & 0                 \\
\hat{\PIPSBl}_1  & 0     & \hdots & \hat{\PIPSBl}_N & 0     & 0                 \\

}},\ K_0^{\text{dense}} := \matr{ \Sigma_0 & A_0\T & \PIPSBlglob_0\T \\ A_0 & 0 & 0 \\ \PIPSBlglob_{0} & 0 & 0 },\
	L^{\text{dense}} := \matr{
	0   & 0      & \PIPSBlglob_1\T   \\
	A_1 & 0      & 0                 \\
	& \vdots &                       \\
	0   & 0      & \PIPSBlglob_N\T   \\
	A_N & 0      & 0                 \\
    \hat{\PIPSBl}_0 & 0     & 0      \\
}
\end{equation*}
in \cref{eq:generalized_sc_system}.
We call the Schur complement system obtained this way the \textit{dense layer} of the multi-layered \HSCA as we do not detect any additional structure in the Schur complement of this layer and thus factorize it with a dense direct linear solver.
Given an implicit factorization of the resulting inner system $K^{\text{dense}}$, we can obtain an implicit factorization of the whole system.
Using a Schur complement to mitigate the fill-in generated by dense columns in \IPMs is usually applied when solving the \textit{normal equations}~\cite{AndersenDenseColumns1996,MeszarosDenseColumns2007}.
When solving the \textit{augmented system}, dense columns are usually treated by adapting the choice of pivoting. The method presented here can be seen as a similar approach~\cite{MarosMeszarosAugmentedSystem1998}, removing both, dense rows and columns from the system.

\subsection{The inner linear systems}\label{sec:inner_linear_systems}

To actually apply \cref{alg:schur_complement_general} to the dense layer, we must be able to solve the inner linear system
\begin{equation}
\label{eq:linear_system_pips_inner_top}
K^{\text{dense}} =
\left[
\begin{array}{ccccc|c}
\Sigma_1        & B_1\T &        &                 &       & \hat{\PIPSBl}_1\T \\
B_1             & 0     &        &                 &       & 0                 \\
                &       & \ddots &                 &       &                   \\
                &       &        & \Sigma_N        & B_N\T & \hat{\PIPSBl}_N\T \\
                &       &        & B_N             & 0     & 0                 \\
\hline
\hat{\PIPSBl}_1 & 0     &        & \hat{\PIPSBl}_N & 0     & 0                 \\
\end{array}
\right]
\matr{
\Delta \primalvector_1 \\
\Delta \dualsEq_1 \\
\vdots \\
\Delta \primalvector_N \\
\Delta \dualsEq_N \\
\Delta \dualsEq_{N+1}
} =
\matr{ r_{\primalvector_1} \\ r_{\dualsEq_1} \\ \vdots \\ r_{\primalvector_N} \\ r_{\dualsEq_N} \\ r_{\dualsEq_{N+1}}}.
\end{equation}
Since \cref{eq:linear_system_pips_inner_top} has the same structure as \cref{eq:linear_system_pips_perturbed_compressed}, we can solve it using the parallel factorization and solve for arrowhead systems introduced in \cref{alg:factor_pips,alg:solve_pips}.
Computing an explicit factorization of the diagonal block matrices would result in a two layered \HSCA system.

To distribute and parallelize the solution process further, we instead split the linear system \cref{eq:linear_system_pips_inner_top} and apply a recursive Schur complement decomposition by permuting subsets of the two-link constraints ``into'' \cref{eq:linear_system_pips_inner_top}.
We recall the structure of \cref{eq:linear_system_pips_inner_top} with the expanded linking constraints:
\begin{equation}
\label{eq:linear_system_pips_inner_top_expanded_links}
\scalemath{0.96}{\left[
\begin{array}{ccccccc|cccc}
\Sigma_1     & B_1\T &               &        &        &                   &       & \PIPSBlloc_1\T    & 0              &        &                       \\
B_1          & 0     &               &        &        &                   &       & 0                 & 0              &        &                       \\
             &       & \Sigma_2      & B_2\T  &        &                   &       & \PIPSBlloc_1'^{T} & \PIPSBlloc_2\T &        &                       \\
             &       & B_2           & 0      &        &                   &       & 0                 & 0              &        &                       \\
             &       &               &        & \ddots &                   &       &                   &                & \ddots &                       \\
             &       &               &        &        & \Sigma_N          & B_N\T &                   &                &        & \PIPSBlloc_{N-1}'^{T} \\
             &       &               &        &        & B_N               & 0     &                   &                &        & 0                     \\
\hline
\PIPSBlloc_1 & 0     & \PIPSBlloc_1' & 0      &        &                   &       &                   &                &        &                       \\
             &       & \PIPSBlloc_2  & 0      &        &                   &       & \mcol{4}{c}{\mrow{3}{*}{0}}                                         \\
             &       &               &        & \ddots &                   &       &                   &                &        &                       \\
             &       &               &        &        & \PIPSBlloc_{N-1}' & 0     &                   &                &        &                       \\
\end{array}
\right] \matr{
\Delta \primalvector_1 \\
\Delta \dualsEq_1 \\
\vdots \\
\Delta \primalvector_N \\
\Delta \dualsEq_N \\
\Delta \dualsEq_{N+1,1} \\
\vdots \\
\Delta \dualsEq_{N+1,N-1}
} =
\matr{ r_{\primalvector_1} \\ r_{\dualsEq_1} \\ \vdots \\ r_{\primalvector_N} \\ r_{\dualsEq_N} \\ r_{\dualsEq_{N+1,1}} \\ \vdots \\ r_{\dualsEq_{N+1,N-1}}}}.
\end{equation}
Next, we group consecutive diagonal blocks of \cref{eq:linear_system_pips_inner_top_expanded_links} into $k\in\N$ subsets by setting $i_1 = 1 < i_2 < \dots < i_{k+1} = N + 1$, $k\in\N$ and assigning to each subset $j \in \{1,\dots,k\}$ the consecutive range of block indices $i_j,\dots,i_{j+1}-1$.
We define
\begin{equation*}
\hat{K}_j :=
\left[
\begin{array}{ccccccc}
	\Sigma_{i_j} & B_{i_j}\T &                &             &        &                     &                \\
    B_{i_j}      & 0         &                &             &        &                     &                \\
                 &           & \Sigma_{i_j+1} & B_{i_j+1}\T &        &                     &                \\
                 &           & B_{{i_j+1}}    & 0           &        &                     &                \\
                 &           &                &             & \ddots &                     &                \\
                 &           &                &             &        & \Sigma_{i_{j+1}-1} & B_{i_{j+1}-1}\T \\
                 &           &                &             &        & B_{i_{j+1}-1}      & 0               \\
\end{array}
\right],
\end{equation*}
\begin{equation*}\scalemath{0.93}{
\Delta \hat{x}_j := \matr{ \Delta \primalvector_{i_j} \\ \Delta\dualsEq_{i_j} \\ \vdots \\ \Delta \primalvector_{i_{j+1}-1} \\ \Delta \dualsEq_{i_{j+1}-1}},\
\Delta \hat{y}_{j} := \matr{
y_{N+1,i_j} \\
\vdots \\
y_{N+1,i_{j+1} - 2} \\
},\
r_{\hat{x}_j} := \matr{ r_{\primalvector_{i_j}} \\ r_{\dualsEq_{i_j}} \\ \vdots \\ r_{\primalvector_{i_{j+1}-1}} \\ r_{\dualsEq_{i_{j+1}-1}}},\
r_{\hat{y}_j} := \matr{ r_{y_{N+1},i_j} \\
\vdots \\
r_{y_{N+1},i_{j+1}-2} \\
}},
\end{equation*}
\begin{equation*}
\hat{\PIPSBlloc}'^{T}_j :=
\left[
\begin{array}{c}
\PIPSBlloc_{i_j}'^{T}\\
0\\
0\\
0\\
\vdots\\
0\\
0\\
\end{array}
\right],\
\hat{B}\T_j :=
\ifthenelse{\zibreport = 1}{\scalemath{0.9}}{}{
\left[
\begin{array}{cccc}
    \PIPSBlloc_{i_j}\T    & 0                    &        &                            \\
    0                     & 0                    &        &                            \\
    \PIPSBlloc_{i_j}'^{T} & \PIPSBlloc_{i_j+1}\T &        &                            \\
                          & 0                    &        &                            \\
                          &                      & \ddots &                            \\
                          &                      &        & \PIPSBlloc_{i_{j+1}-2}'^{T}\\
                          &                      &        & 0                          \\
\end{array}
\right]},\
\hat{\PIPSBlloc}_j^{T} :=
\left[
\begin{array}{c}
	0\\
	0\\
	0\\
	0\\
	\vdots\\
	\PIPSBlloc_{i_{j+1}-1}\T\\
	0\\
\end{array}
\right],
\end{equation*}
for $j=1,\dots,k$ ($\hat{\PIPSBlloc}'^{T}_j$, $\hat{\PIPSBlloc}_j^{T}$ are only defined up to $k-1$). We can regroup \cref{eq:linear_system_pips_inner_top_expanded_links} into
\begin{equation}
\label{eq:inner_system_before_reorder}
\scalemath{0.99}{
\left[
\begin{array}{cccc|cccccc}
    \hat{K}_1             &                         &        &                           & \hat{B}\T_1 & \hat{\PIPSBlloc}_{1}\T    &             &        &                               &            \\
                          & \hat{K}_2               &        &                           &             & \hat{\PIPSBlloc}_{1}^{T'} & \hat{B}\T_2 & \hat{\PIPSBlloc}_{2}\T       &                               &            \\
                          &                         & \ddots &                           &             &                           &             & \mcol{2}{c}{\ddots}                               &            \\
                          &                         &        & \hat{K}_k                 &             &                           &             &        & \hat{\PIPSBlloc}_{k-1}^{T'} & \hat{B}\T_k\\
    \hline
    \hat{B}_1             &                         &        &                           &             &                           &             &        &                               &            \\
    \hat{\PIPSBlloc}_{1}  & \hat{\PIPSBlloc}_{1}'   &        &                           &             &                           &             &        &                               &            \\
                          & \hat{B}_2               &        &                           & \mcol{6}{c}{\mrow{3}{*}{0}} \\
                          & \hat{\PIPSBlloc}_{2}    & \mrow{2}{*}{\ddots} &                           &             &                           &             &        &                               &            \\
                          &                         &        & \hat{\PIPSBlloc}_{k-1}'   &             &                           &             &        &                               &            \\
                          &                         &        & \hat{B}_k                 &             &                           &             &        &                               &            \\
\end{array}
\right] \matr{
\Delta \hat{x}_{1} \\
\Delta \hat{x}_2 \\
\vdots \\
\Delta \hat{x}_k \\
\Delta \hat{y}_{1} \\
\Delta \dualsEq_{N+1, i_{2}-1} \\
\Delta \hat{y}_{2} \\
\vdots \\
\Delta \dualsEq_{N+1, i_{k}-1} \\
\Delta \hat{y}_{k} \\
} = \matr{
r_{\hat{x}_1} \\
r_{\hat{x}_2} \\
\vdots \\
r_{\hat{x}_k} \\
r_{\hat{y}_1} \\
r_{y_{N+1}, i_{2}-1} \\
r_{\hat{y}_2} \\
\vdots \\
r_{y_{N+1}, i_{k}-1} \\
r_{\hat{y}_k}\\
}}.
\end{equation}
A symmetric permutation of \cref{eq:inner_system_before_reorder} gives
\begin{equation}\label{eq:inner_system_after_reorder}
\scalemath{0.99}{
\left[
\begin{array}{ccccccc|ccc}
    \hat{K}_1             & \hat{B}\T_1 &                         &             &        &                           &             & \hat{\PIPSBlloc}_{1}\T        &        &                               \\
    \hat{B}_1             & 0           &                         &             &        &                           &             & 0                             &        &                               \\
                          &             & \hat{K}_2               & \hat{B}\T_2 &        &                           &             & \hat{\PIPSBlloc}_{1}^{T'}     & \hat{\PIPSBlloc}_{2}\T       &                               \\
                          &             & \hat{B}_2               & 0           &       &                           &             & 0                             & 0       &                               \\
                          &             &                         &             & \ddots &                           &             &                               & \mcol{2}{c}{\ddots}                               \\
                          &             &                         &             &        & \hat{K}_k                 & \hat{B}\T_k &                               &        & \hat{\PIPSBlloc}_{k-1}^{T'}   \\
                          &             &                         &             &        & \hat{B}_k                 & 0           &                               &        & 0                             \\
    \hline
    \hat{\PIPSBlloc}_{1}  &0            & \hat{\PIPSBlloc}_{1}'   & 0           &        &                           &             &                               &        &                               \\
                          &             &  \hat{\PIPSBlloc}_{2}      & 0           & \mrow{2}{*}{\ddots} &                           &             & \mcol{3}{c}{\mrow{1}{*}{0}} \\
                          &             &                         &             &        & \hat{\PIPSBlloc}_{k-1}'   & 0           &                               &        &                               \\
\end{array}
\right] } \matr{
\Delta \hat{x}_{1} \\
\Delta \hat{y}_{1} \\
\vdots \\
\Delta \hat{x}_k \\
\Delta \hat{y}_{k} \\
\Delta \dualsEq_{N+1, i_{2}-1} \\
\vdots \\
\Delta \dualsEq_{N+1, i_{k}-1} \\
} = \matr{
r_{\hat{x}_1} \\
r_{\hat{y}_1} \\
\vdots \\
r_{\hat{x}_k} \\
r_{\hat{y}_k} \\
r_{y_{N+1}, i_{2}-1} \\
\vdots \\
r_{y_{N+1}, i_{k}-1} \\
}.
\end{equation}
By setting
\begin{equation*}
K_i^{inner} := \matr{
    \hat{K}_i & \hat{B}_i\T \\
    \hat{B}_i & 0     \\
},\quad K^{inner}_0 := 0,\quad L_i^{inner} := \matr{
    \dots & 0 & \hat{\PIPSBlloc}_{i-1}^{T'} & \hat{\PIPSBlloc}_{i}\T & 0 & \dots \\
    \dots     & 0 & 0                           & 0                    & 0 & \dots \\
},
\end{equation*}
the system in \cref{eq:inner_system_after_reorder} is of the same form as \cref{eq:linear_system_pips_perturbed_compressed}, and we can use \cref{alg:factor_pips} and \cref{alg:solve_pips} to factorize and solve it in parallel.
However, instead of computing a direct factorization of each diagonal $2\times2$-block, we notice their arrowhead structure:
\begin{equation}\label{eq:linear_system_pips_inner_lower}
K_i^{inner} = \scalemath{1.0}{\left[\begin{array}{ccccccc|cccc}
\Sigma_{i_j} & B_{i_j}\T &                &             &        &                     &                & \PIPSBlloc_{i_j}\T    & 0                    &        &                            \\
B_{i_j}      & 0         &                &             &        &                     &                & 0                     & 0                    &        &                            \\
             &           & \Sigma_{i_j+1} & B_{i_j+1}\T &        &                     &                & \PIPSBlloc_{i_j}^{T'} & \PIPSBlloc_{i_j+1}\T &        &                            \\
             &           & B_{{i_j+1}}    & 0           &        &                     &                & 0                     & 0                    &        &                            \\
             &           &                &             & \ddots &                     &                &                       &                      & \ddots &                            \\
             &           &                &             &        & \Sigma_{i_{j+1}-1} & B_{i_{j+1}-1}\T &                       &                      &        & \PIPSBlloc_{i_{j+1}-2}^{T'}\\
             &           &                &             &        & B_{i_{j+1}-1}      & 0               &                       &                      &        & 0                          \\
\hline
\PIPSBlloc_{i_j} & 0 & \PIPSBlloc_{i_j}^{'}& 0         &        &                     &                &                       &                      &        &                            \\
0            & 0 & \PIPSBlloc_{i_j+1}   & 0           &        &                     &                &  \mcol{4}{c}{\mrow{3}{*}{0}} \\
             &           &                &             & \ddots &                     &                &                       &                      &        &                            \\
             &           &                &             &        & \PIPSBlloc_{i_{j+1}-2}^{'}& 0       &                       &                      &        &                            \\
\end{array}\right]},
\end{equation}
$i=1,\dots,k$.
Each diagonal block itself is of the form \cref{eq:linear_system_pips_inner_top_expanded_links} and can again be handled by algorithms \cref{alg:factor_pips} and \cref{alg:solve_pips} to obtain an implicit factorization of $K^{inner}_i$.
Alternatively, we can split each $K^{inner}_i$ further by applying the same technique as outlined above, recursively.

We have decomposed the linear system \cref{eq:linear_system_pips_inner_top} into $k$ smaller linear systems defined by the system matrices $K^{inner}_i$.
We illustrate this process schematically in \cref{fig:hierarchical_decomposition}.
\begin{figure}
    \centering
    \begin{subfigure}{.35\textwidth}
    \centering
    \begin{tikzpicture}[scale=0.65]
\draw[cp2, fill=cp2!50] (0,0) rectangle +(0.5,0.5);
\draw[cp2, fill=cp2!50] (-0.5,0.5) rectangle +(0.5,0.5);
\draw[cp2, fill=cp2!50] (-1,1) rectangle +(0.5,0.5);
\draw[cp5, fill=cp5!50] (-1.5,1.5) rectangle +(0.5,0.5);
\draw[cp5, fill=cp5!50] (-2,2) rectangle +(0.5,0.5);
\draw[cp5, fill=cp5!50] (-2.5,2.5) rectangle +(0.5,0.5);
\draw[cp1, fill=cp1!70] (-3,3) rectangle +(0.5,0.5);
\draw[cp1, fill=cp1!70] (-3.5,3.5) rectangle +(0.5,0.5);
\draw[cp1, fill=cp1!70] (-4,4) rectangle +(0.5,0.5);

\draw[cp2, fill=cp2!50] (-0.5,-0.8) rectangle +(1,0.1);
\draw[cp2, fill=cp2!50] (-1,-0.7) rectangle +(1,0.1);
\draw[cp3, fill=cp3!80] (-1.5,-0.6) rectangle +(1,0.1);
\draw[cp5, fill=cp5!50] (-2,-0.5) rectangle +(1,0.1);
\draw[cp5, fill=cp5!50] (-2.5,-0.4) rectangle +(1,0.1);
\draw[cp3, fill=cp3!80] (-3,-0.3) rectangle +(1,0.1);
\draw[cp1, fill=cp1!70] (-3.5,-0.2) rectangle +(1,0.1);
\draw[cp1, fill=cp1!70] (-4,-0.1) rectangle +(1,0.1);

\draw[cp2, fill=cp2!50] (1.2,0.0) rectangle +(0.1,1);
\draw[cp2, fill=cp2!50] (1.1,0.5) rectangle +(0.1,1);
\draw[cp3, fill=cp3!80] (1.0,1.0) rectangle +(0.1,1);
\draw[cp5, fill=cp5!50] (0.9,1.5) rectangle +(0.1,1);
\draw[cp5, fill=cp5!50] (0.8,2.0) rectangle +(0.1,1);
\draw[cp3, fill=cp3!80] (0.7,2.5) rectangle +(0.1,1);
\draw[cp1, fill=cp1!70] (0.6,3.0) rectangle +(0.1,1);
\draw[cp1, fill=cp1!70] (0.5,3.5) rectangle +(0.1,1);

\draw[black, dashed] (-3, 0.0) -- (1.2, 0.0); 
\draw[black, dashed] (0.5, 3.5) -- (0.5, -0.7); 

\draw[black] (-4.0,4.5) -- (-4.0,-0.8); 
\draw[black] (1.3,4.5) -- (1.3,-0.8); 
\draw[black] (-4.0,4.5) -- (1.3,4.5); 
\draw[black] (-4.0, -0.8) -- (1.3, -0.8); 

\end{tikzpicture}
    \caption{Inner system \cref{eq:linear_system_pips_inner_top_expanded_links}}
    \label{fig:inner_system_before_split}
    \end{subfigure}%
    \begin{subfigure}{.35\textwidth}
    \centering
    \begin{tikzpicture}[scale=0.65]
\draw[cp2, fill=cp2!50] (0,0) rectangle +(0.5,0.5);
\draw[cp2, fill=cp2!50] (-0.5,0.5) rectangle +(0.5,0.5);
\draw[cp2, fill=cp2!50] (-1,1) rectangle +(0.5,0.5);
\draw[cp2, fill=cp2!50] (-1,-0.1) rectangle +(1,0.1);
\draw[cp2, fill=cp2!50] (-0.5,-0.2) rectangle +(1,0.1);
\draw[cp2, fill=cp2!50] (0.5,0.5) rectangle +(0.1,1);
\draw[cp2, fill=cp2!50] (0.6,0.0) rectangle +(0.1,1);

\draw[cp5, fill=cp5!50] (-1.7,1.7) rectangle +(0.5,0.5);
\draw[cp5, fill=cp5!50] (-2.2,2.2) rectangle +(0.5,0.5);
\draw[cp5, fill=cp5!50] (-2.7,2.7) rectangle +(0.5,0.5);
\draw[cp5, fill=cp5!50] (-2.7,1.6) rectangle +(1,0.1);
\draw[cp5, fill=cp5!50] (-2.2,1.5) rectangle +(1,0.1);
\draw[cp5, fill=cp5!50] (-1.2,2.2) rectangle +(0.1,1);
\draw[cp5, fill=cp5!50] (-1.1,1.7) rectangle +(0.1,1);

\draw[cp1, fill=cp1!70] (-3.4,3.4) rectangle +(0.5,0.5);
\draw[cp1, fill=cp1!70] (-3.9,3.9) rectangle +(0.5,0.5);
\draw[cp1, fill=cp1!70] (-4.4,4.4) rectangle +(0.5,0.5);
\draw[cp1, fill=cp1!70] (-4.4,3.3) rectangle +(1,0.1);
\draw[cp1, fill=cp1!70] (-3.9,3.2) rectangle +(1,0.1);
\draw[cp1, fill=cp1!70] (-2.9,3.9) rectangle +(0.1,1);
\draw[cp1, fill=cp1!70] (-2.8,3.4) rectangle +(0.1,1);

\draw[cp3, fill=cp3!50] (-3.2,-0.3) rectangle +(1,0.1);
\draw[cp3, fill=cp3!50] (-1.5,-0.4) rectangle +(1,0.1);

\draw[cp3, fill=cp3!50] (0.7,2.7) rectangle +(0.1,1);
\draw[cp3, fill=cp3!50] (0.8,1.0) rectangle +(0.1,1);

\draw[black] (-4.4,4.9) -- (-4.4,-0.4); 
\draw[black] (0.9,4.9) -- (0.9,-0.4); 
\draw[black] (-4.4,4.9) -- (0.9,4.9); 
\draw[black] (-4.4, -0.4) -- (0.9, -0.4); 

\end{tikzpicture}
    \caption{Split inner system \cref{eq:inner_system_after_reorder}}
    \label{fig:inner_system_split}
    \end{subfigure}%
    \begin{subfigure}{.35\textwidth}
    \centering
    \begin{tikzpicture}[scale=0.65]
\draw[cp2, fill=cp2!50] (-1,-0.2) rectangle +(1.7,1.7);
\draw[cp5, fill=cp5!50] (-2.7,1.5) rectangle +(1.7,1.7);
\draw[cp1, fill=cp1!70] (-4.4,3.2) rectangle +(1.7,1.7);

\draw[cp3, fill=cp3!80] (-3.2,-0.3) rectangle +(1,0.1);
\draw[cp3, fill=cp3!80] (-1.5,-0.4) rectangle +(1,0.1);

\draw[cp3, fill=cp3!80] (0.7,2.7) rectangle +(0.1,1);
\draw[cp3, fill=cp3!80] (0.8,1.0) rectangle +(0.1,1);

\draw[black] (-4.4,4.9) -- (-4.4,-0.4); 
\draw[black] (0.9,4.9) -- (0.9,-0.4); 
\draw[black] (-4.4,4.9) -- (0.9,4.9); 
\draw[black] (-4.4, -0.4) -- (0.9, -0.4); 

\end{tikzpicture}
    \caption{First inner layer.}
    \label{fig:first_inner_layer}
    \end{subfigure}%
    \caption{Permuting the inner linear system}
    \label{fig:hierarchical_decomposition}
\end{figure}
In \cref{fig:inner_system_before_split}, we are given a system $K^{dense}$ with nine diagonal blocks and expanded two-link structure.
We sub-divide this system into three smaller ones.
For this purpose, we color-coded the components belonging to each subsystem with the same color.
By permuting all local linking constraints but the ones linking the subsystems (depicted in green) further into the system, we obtain the system in \cref{fig:inner_system_split}, with arrowhead systems ${K}^{inner}_i$, $i=1,2,3$, on the diagonal.
Interpreting the system in \cref{fig:inner_system_split} as displayed in \cref{fig:first_inner_layer}, we can use \cref{alg:factor_pips} for its factorization, assuming an implicit factorization for ${K}^{inner}_i$ is available.
When going from \cref{eq:inner_system_before_reorder} to \cref{eq:inner_system_after_reorder} it is not necessary that all linking constraints permuted \emph{into} the system are two-links.
For $\hat{\PIPSBlloc}_j$ and $\hat{\PIPSBlloc}'_j$, $j=1,\dots,k-1$, we can allow linking constraints connecting at most the diagonal blocks contained in the range $i_j,\dots,i_{j+1}-1$.
%

Recursively applying the splitting procedure to each ${K}_i^{inner}$ would create another layer in the Schur complement decomposition and further shrink the lowest-level systems.
Given enough blocks, this enables an arbitrary number of layers.
In practice, the recursive Schur complement decomposition affects the computational cost of a linear system solve and must be weighed against the benefit of smaller Schur complements and their distributed computation.
We have not encountered systems or problems where more than four layers, one dense layer and three inner layers, were beneficial.

Note that the distribution of the linear systems and their subsystems naturally groups the \MPI processes into subsets defined via the blocks assigned to the systems ${K}_i^{inner}$.
\begin{figure}[htbp]\label{fix:mpi_assignment}
\centering

\tikzset{
  node/.style={draw=none, fill=none},
  proc/.style={rectangle, draw=none, text=cp2},
  wrong/.style={rectangle, draw=none, text=cp6},
  dashedarrow/.style={dashed, ->, thick}
}

\begin{subfigure}[b]{0.48\textwidth}
\centering
\begin{tikzpicture}[
  scale=0.9, transform shape,
  level distance=0.8cm,
  level 1/.style={sibling distance=3cm},
  level 2/.style={sibling distance=1.1cm}
]
\node[node] (R) {$K^{dense}$}
  child { node[node] (A) {${K}^{inner}_1$}
    child { node[node] (B0) {$B_1$} }
    child { node[node] (B1) {$B_2$} }
    child { node[node] (B2) {$B_3$} }
  }
  child { node[node] (B) {${K}^{inner}_2$}
    child { node[node] (B3) {$B_4$} }
    child { node[node] (B4) {$B_5$} }
    child { node[node] (B5) {$B_6$} }
  };

\node[proc] at ($(B0) + (0,-0.6)$) (P0) {$P_0$};
\node[proc] at ($(B1) + (0,-0.6)$) (P1) {$P_1$};
\node[proc] at ($(B2) + (0,-0.6)$) (P2) {$P_2$};
\node[proc] at ($(B3) + (0,-0.6)$) (P3) {$P_3$};
\node[proc] at ($(B4) + (0,-0.6)$) (P4) {$P_4$};
\node[proc] at ($(B5) + (0,-0.6)$) (P5) {$P_5$};
\end{tikzpicture}
\caption{One process per block (\color{cp2}valid).}
\label{fig:mpi_assignment_1to1}
\end{subfigure}%
\begin{subfigure}[b]{0.48\textwidth}
\centering
\begin{tikzpicture}[
  scale=0.9, transform shape,
  level distance=0.8cm,
  level 1/.style={sibling distance=3cm},
  level 2/.style={sibling distance=1.1cm}
]
\node[node] (R) {$K^{dense}$}
  child { node[node] (A) {${K}^{inner}_1$}
    child { node[node] (B0) {$B_1$} }
    child { node[node] (B1) {$B_2$} }
    child { node[node] (B2) {$B_3$} }
  }
  child { node[node] (B) {${K}^{inner}_2$}
    child { node[node] (B3) {$B_4$} }
    child { node[node] (B4) {$B_5$} }
    child { node[node] (B5) {$B_6$} }
  };

\node[wrong] at ($(B0) + (0,-0.6)$) (P0) {$P_0$};
\node[wrong] at ($(B1) + (0,-0.6)$) (P0b) {$P_0$};
\node[wrong] at ($(B2) + (0,-0.6)$) (P1) {$P_1$};
\node[wrong] at ($(B3) + (0,-0.6)$) (P1b) {$P_1$};
\node[wrong] at ($(B4) + (0,-0.6)$) (P2) {$P_2$};
\node[wrong] at ($(B5) + (0,-0.6)$) (P2b) {$P_2$};

\end{tikzpicture}
\caption{$P_1$ spans two subsystems ({\color{cp6}invalid}).}
\label{fig:mpi_assignment_multiple_shared_trees}
\end{subfigure}
\begin{subfigure}[b]{0.48\textwidth}
\centering
\begin{tikzpicture}[
  scale=0.9, transform shape,
  level distance=0.8cm,
  level 1/.style={sibling distance=3cm},
  level 2/.style={sibling distance=1.1cm}
]
\node[node] (R) {$K^{dense}$}
  child { node[node] (A) {${K}^{inner}_1$}
    child { node[node] (B0) {$B_1$} }
    child { node[node] (B1) {$B_2$} }
    child { node[node] (B2) {$B_3$} }
  }
  child { node[node] (B) {${K}^{inner}_2$}
    child { node[node] (B3) {$B_4$} }
    child { node[node] (B4) {$B_5$} }
    child { node[node] (B5) {$B_6$} }
  };

\node[proc] at ($(B0) + (0,-0.6)$) (P0) {$P_0$};
\node[proc] at ($(B1) + (0,-0.6)$) (P0b) {$P_0$};
\node[proc] at ($(B2) + (0,-0.6)$) (P0c) {$P_0$};
\node[proc] at ($(B3) + (0,-0.6)$) (P1) {$P_1$};
\node[proc] at ($(B4) + (0,-0.6)$) (P1b) {$P_1$};
\node[proc] at ($(B5) + (0,-0.6)$) (P1c) {$P_1$};
\end{tikzpicture}
\caption{One process per subsystem (\color{cp2}valid).}
\label{fig:mpi_assignment_fully_owned_subsystem}
\end{subfigure}%
\begin{subfigure}[b]{0.48\textwidth}
\centering
\begin{tikzpicture}[
  scale=0.9, transform shape,
  level distance=0.8cm,
  level 1/.style={sibling distance=3cm},
  level 2/.style={sibling distance=1.1cm}
]
\node[node] (R) {$K^{dense}$}
  child { node[node] (A) {${K}^{inner}_1$}
    child { node[node] (B0) {$B_1$} }
    child { node[node] (B1) {$B_2$} }
    child { node[node] (B2) {$B_3$} }
  }
  child { node[node] (B) {${K}^{inner}_2$}
    child { node[node] (B3) {$B_4$} }
    child { node[node] (B4) {$B_5$} }
    child { node[node] (B5) {$B_6$} }
  };

\node[proc] at ($(B0) + (0,-0.6)$) (P0) {$P_0$};
\node[proc] at ($(B1) + (0,-0.6)$) (P1b) {$P_0$};
\node[proc] at ($(B2) + (0,-0.6)$) (P1) {$P_1$};
\node[proc] at ($(B3) + (0,-0.6)$) (P2) {$P_2$};
\node[proc] at ($(B4) + (0,-0.6)$) (P2b) {$P_2$};
\node[proc] at ($(B5) + (0,-0.6)$) (P3) {$P_3$};
\end{tikzpicture}
\caption{Processes share single subsystem (\color{cp2}valid).}
\label{fig:mpi_assignment_single_shared_subsystem}
\end{subfigure}

\caption{System factorization/solve tree: {\color{cp2}valid} and {\color{cp6}invalid} process assignments.}
\label{fig:mpi-assignment}
\end{figure}
We schematically display the assignment of different numbers of processes to different subsystems in \cref{fig:mpi-assignment}.
We display the factorization/solve tree of a two level \HSCA system, splitting the inner system into two subsystems and omitting the dense layer.
We denote the problem blocks representatively by $B_i$, $i=1,\dots,6$.
Assignment of a process $i=0,\dots,5$ (\MPI assigns processes ranks starting from 0) to a block is visualized by placing $P_i$ below the respective block in the tree.
In \cref{fig:mpi_assignment_1to1} processes $P_0$, $P_1$, and $P_2$ are assigned to ${K}^{inner}_1$ and $P_3$, $P_4$, and $P_5$ to ${K}^{inner}_2$.
When splitting linear systems into smaller subsystems, the processes associated with the subsystem's blocks are grouped together and communally compute the implicit factorization of the subsystem.
Within \PIPS, this grouping is realized using \MPI groups and communicators, which simplifies the setup and the data transfer required for subsystem factorization and solve.
This assignment becomes non-trivial when fewer than $N$ processes are available.
The parallelization of \HSCA relies on the subsystems being factorized and used for solving linear equations in parallel.
If a single process is responsible for factorizing two linear systems in the same level, it will have to communicate with the other processes responsible for each of the systems.
This situation is depicted in \cref{fig:mpi_assignment_multiple_shared_trees}, where $P_1$ is assigned to ${K}^{inner}_1$ and ${K}^{inner}_2$.
As \MPI communication operations are blocking until all processes are available, this assignment blocks the parallel factorization of ${K}^{inner}_1$ and ${K}^{inner}_2$ and, at least partially, renders the computations sequential.
Thus, an assignment of processes to blocks must consider the splitting within \HSCA.
As a result, in the current development, we have implemented a simple heuristic that enforces two rules: either one process is fully responsible for any linear system of a given level it is assigned to, as in \cref{fig:mpi_assignment_fully_owned_subsystem}, or it participates in the processing of exactly one subsystem as in \cref{fig:mpi_assignment_single_shared_subsystem}.
\section{Computational experiments}
\label{sec:comp_exps}

In this section, we describe numerical results obtained with \HSCA.
All experiments with \PIPS were conducted on the \JUWELS supercomputer \cite{JUWELS} at Forschungszentrum J{\"{u}}lich.
\JUWELS has 2,271 standard compute nodes (Dual Intel Xeon Platinum 8168), each with 96GB memory and 2x24 cores.
The nodes are connected via a Mellanox EDR InfiniBand high-speed network.
The experiments with other commercial solvers were conducted on a compute cluster (Intel Xeon Gold 6338) at the Zuse Institute Berlin, each with 1024GB memory and 2x32 cores.
As this paper is not a commercial solver comparison, we have anonymized all results obtained with commercial optimization software.
Our software is freely available on \href{https://github.com/NCKempke/PIPS-IPMpp}{GitHub}.

\PIPS is able to not only leverage parallelism via \MPI but also shared memory parallelism via \OpenMP \cite{dagum1998_openmp}.
\OpenMP is mainly used for the factorization and the solution of the linear Schur complement systems within the parallel solver \PARDISO (see \cite{Rehfeldt2019_PIPSIPMpp}).
In preliminary experiments with \HSCA, we compared the performance of the linear solvers \MAEightySix \cite{MA86Hogg2010AnIS}, \MAFiftySeven \cite{MA57}, \MATwentySeven \cite{MA27}, \MUMPS \cite{MUMPS:1}, \PARDISO \cite{pardiso-7.2a}, \MKL \PARDISO \cite{MKL_PARDISO}, and \WSMP \cite{WSMP}.
Our initial studies focused on two key aspects: the speed of factorization for medium-size linear systems and strong scalability when solving a given linear system with multiple right-hand sides, which is the most computationally demanding part of \HSCA.
We finally settled on \MAFiftySeven.
This deviates from the latest incomplete LU factorization approach applied in \cite{Petra2014_PIPSaugmented,Rehfeldt2019_PIPSIPMpp} and \PIPS' original implementation.
While the factorization routine of \MAFiftySeven is slower than that of \PARDISO, particularly for larger systems, it allows for more efficient parallel solves.
Specifically, the sequential \MAFiftySeven is thread-safe, enabling the use of \OpenMP to solve multiple right-hand sides simultaneously by calling the solver from different threads using the same factorization.
In our experience, this ``outside'' multi-threading of the solve routine scales much better than the solve routines of all other solvers we tested.
Additionally, by not relying on\PARDISO's parallelized incomplete LU factorization for computing the bottommost Schur complement contributions, we gain numerical precision on the lowest \HSCA level, albeit at the cost of losing some of \PARDISO's speed-up.
\PARDISO on the other hand showed greater robustness when solving rank-deficient linear systems.
In summary, for medium-size matrices where repeated linear solves dominate the computation, the slower factorization of \MAFiftySeven is outweighed by its superior scalability and precision in solving multiple right-hand sides. These advantages made \MAFiftySeven the preferred solver for our setting.

\subsection*{The scaling behavior of the hierarchical approach}

We first assess the scaling behavior of our implementation.
In \cref{fig:miso_scaling} we show the behavior of our solver on the \textit{MISO\_DISP\_488} instance as well as the scaling behavior of the three commercial but academically available solvers \CPLEX 22.1.1.0 \cite{CPLEX12}, \GUROBI 11 \cite{Gurobi11}, and the Cardinal Optimizer \COPT 7.2.3 \cite{COPT71}.
\begin{figure}
\centering
\scalebox{.9}{\begin{tikzpicture}
    \begin{loglogaxis}
    [
        xlabel={\#Threads (for Solvers)/\#MPI processes},
        ylabel={Time to solution [s]},
        log ticks with fixed point,
        scale only axis,
        xmin=8, xmax=4000,
        ymin=100, ymax=10000,
        width=11cm,
        height=7cm,
        legend pos = south west,
        grid style={line width=1pt, draw=black, opacity=0.1},
        xmajorgrids={true},
        ymajorgrids={true},
        log basis x=2
    ]
    \addplot
    [
        color=cp6,
        mark=square,
        thick
    ]
    coordinates {
    (92, 5055)(184, 2544)(368, 1385)(736,735)(1472,482)(2190,379)
    };
    \addplot
    [
        color=cp2,
        mark=square*,
        thick
    ]
    coordinates {
    (92,2785)(184, 1341)(368, 739)(736,404)(1472,250)(2190,200)
    };
    \addplot
    [
        color=purple,
        mark=triangle,
        thick
    ]
    coordinates {
    (1, 2468)(2,2958)(4,2304)(8,3098)(16,1527)(32,1373)(64,1866)
    };
    \addplot
    [
        color=cp5,
        mark=diamond,
        thick
    ]
    coordinates {
    (1, 18648)(2,12297)(4,7131)(8,3602)(16,3542)(32,3363)(64,3962)
    };
    \addplot
    [
        color=cp4,
        mark=pentagon,
        thick
    ]
    coordinates {
    (1, 2702)(2,2414)(4,1618)(8,1259)(16,1092)(32,1060)(64,1109)
    };
    \addplot
    [
        color=cp3,
        mark=*,
        style=dashed,
        ultra thin
    ]
    coordinates {
    (92, 5055)(184, 2527)(368, 1263)(736,632)(1472,316)(2190,212)
    };
    \addplot
    [
        color=cp3,
        mark=*,
        style=dashed,
        ultra thin
    ]
    coordinates {
    (92,2785)(184, 1393)(368, 696)(736,348)(1472,174)(2190,117)
    };
    \addlegendentry{1 \OpenMP thread}
    \addlegendentry{2 \OpenMP threads}
    \addlegendentry{Solver 1}
    \addlegendentry{Solver 2}
    \addlegendentry{Solver 3}
    \addlegendentry{Ideal speed-up}
    \end{loglogaxis}
\end{tikzpicture}}
\caption{Scaling behavior of all solvers on MISO\_DISP\_488.}
\label{fig:miso_scaling}
\end{figure}
We solved the instance with \PIPS multiple times using a different number of \MPI processes and \OpenMP threads and (implicitly) cluster nodes.
The problem displayed has 120 million nonzeros, 25 million columns, and 12 million rows and is split into 2,190 blocks, which limits the number of usable \MPI processes to 2,190.
Notably, the instance contains roughly 55,000 linking constraints, and the initial symbolic analysis of the Schur complement in \PIPS determined an upper bound of 100,424,258 nonzeros.
This made the instance intractable for the original \PIPS algorithm (see \cref{tab:results}).
We adapted the number of threads for each solution with the commercial solvers while turning off crossover and using barrier as the optimization algorithm.
Both \PIPS and the commercial solvers were set to solve the \LP up to a convergence tolerance of $10^{-6}$.

Additionally, we plot the optimal speed-up for each scenario as a dotted line.
Our approach is not expected to achieve linear speed-up, as the mid-level Schur complements are solved by only a subset of threads in parallel, and the root-level Schur complement is solved sequentially.
However, \HSCA achieves a strong near-linear speed-up on the given instance.
Additionally, for this instance \PIPS benefits from adding a second \OpenMP thread while solving.
More \OpenMP threads did not produce much speed-up over the 2-thread variant and the usefulness of \OpenMP is generally problem dependent.
It can also be observed, that the commercial solvers do not scale well with the availability of more compute resources.
Each solver runs fastest with 32 threads, which we attribute to the memory-boundedness of the Cholesky factorization underlying the \IPM.

\subsection*{Impact of block sizes and splitting on the performance of hierarchical approach vs. original approach}

While we have previously demonstrated near-ideal scaling on instances suited for parallelization, in this section we deliberately examine instances where the growth of the Schur complement dominates.
This highlights the strong-scaling trade-offs and limitations inherent in both the original and hierarchical approaches.
We display the different negative impact to be expected from growing Schur complements for \PIPS and \HSCA.

The performance of both the original \PIPS and \HSCA is heavily problem-dependent.
Given the superior performance of \PARDISO on large linear systems, \PIPS is largely favored on systems with a manageable number of linking constraints and large blocks of constraints.
On the other hand, the hierarchical extension struggles with large Schur complements and large diagonal blocks since \MAFiftySeven's factorization routine is slower on large systems.

Additionally, on a given problem the performance of either approach is decomposition dependent.
A ``good'' decomposition balances two competing effects: coarser blocks reduce the size of the Schur complement but increase diagonal block factorization time, whereas finer blocks decrease diagonal block factorization time but increase both the size of the global Schur complement and communication overhead among MPI processes.
For the original approach, the time spent in the Schur complement is highly sensitive to a growing number of linking constraints, whereas for \HSCA, this impact is mitigated by the additional Schur complement parallelization.
We note that for neither approach does using the finest possible splitting yield the best runtime.
While finer splitting reduces the size of diagonal blocks, it enlarges the global Schur complement and increases communication overhead among the growing number of MPI processes. As a result, runtime can increase with more processors---an expected strong-scaling limitation in Schur complement methods.

In our view, the two approaches complement one another, each being efficient for distinct problem classes.
The performance will depend on how a modeler or automatic structure detection algorithm decomposes a given model into blocks: the original approach will favor few large diagonal blocks and a manageable Schur complement; \HSCA will favor more and smaller diagonal blocks even if this results in larger Schur complements.

To demonstrate this, we conducted a simple experiment on four small to medium \SIMPLE instances with different sizes (\cref{fig:splitting}).
These four instances illustrate the behavior under Schur complement stress.
From the top left to the bottom right, the instances become increasingly more difficult to solve, simply due to a growing size of their diagonal blocks.
The number of linking constraints in the problem, when split into a certain number of blocks, remains constant over all instances.
Each instance can theoretically be split into 8,760 blocks.
Combining multiple blocks into one decreases the number of linking constraints in the Schur complement approach.
We solved each instance with the original \PIPS and our new approach, split into five different amounts of blocks: 365, 730, 1,460, 2,920, and 8,760.
Here, we always use as many \MPI processes as blocks available.
In this setting, a lower number of blocks and processes mean bigger diagonal blocks and, as pointed out before, smaller global Schur complements.

\begin{figure}
\begin{subfigure}{.5\textwidth}
\caption{nb5\_to1}
\begin{adjustbox}{width=\textwidth}
\begin{tikzpicture}
    \begin{loglogaxis}
    [
        xlabel={\#blocks/\MPI processes},
        ylabel={Time to solution [s]},
        log ticks with fixed point,
        scale only axis,
        xmin=300, xmax=10000,
        ymin=1, ymax=1000,
        legend pos=south east,
        grid style={line width=1pt, draw=black, opacity=0.1},
        xmajorgrids={true},
        ymajorgrids={true},
        log basis x=2
    ]
    \addplot
    [
        color=cp2,
        mark=square,
        thick
    ]
    coordinates {
    (365,3.5)(730,7)(1460,17)(2920,36)(8760,174)
    };
    \addplot
    [
        color=cp6,
        mark=square*,
        thick
    ]
    coordinates {
    (365,2.4)(730,3)(1460,4.8)(2920,5.1)(8760,15.4)
    };
    \addlegendentry{\PIPS}
    \addlegendentry{Hierarchical}
    \end{loglogaxis}
\end{tikzpicture}
\end{adjustbox}
\end{subfigure}%
\begin{subfigure}{.5\textwidth}
\caption{nb15\_to1}
\begin{adjustbox}{width=\textwidth}
\begin{tikzpicture}
    \begin{loglogaxis}
    [
        xlabel={\#blocks/\MPI processes},
        ylabel={Time to solution [s]},
        log ticks with fixed point,
        scale only axis,
        xmin=300, xmax=10000,
        ymin=1, ymax=1000,
        legend pos=south east,
        grid style={line width=1pt, draw=black, opacity=0.1},
        xmajorgrids={true},
        ymajorgrids={true},
        log basis x=2
    ]
    \addplot
    [
        color=cp2,
        mark=square,
        thick
    ]
    coordinates {
    (365,14)(730,21)(1460,40)(2920,86)(8760,320)
    };
    \addplot
    [
        color=cp6,
        mark=square*,
        thick
    ]
    coordinates {
    (365,13.5)(730,12.6)(1460,12.4)(2920,17.8)(8760,39.7)
    };
    \addlegendentry{\PIPS}
    \addlegendentry{Hierarchical}
    \end{loglogaxis}
\end{tikzpicture}
\end{adjustbox}
\end{subfigure}
\begin{subfigure}{.5\textwidth}
\caption{nb30\_to1}
\begin{adjustbox}{width=\textwidth}
\begin{tikzpicture}
    \begin{loglogaxis}
    [
        xlabel={\#blocks/\MPI processes},
        ylabel={Time to solution [s]},
        log ticks with fixed point,
        scale only axis,
        xmin=300, xmax=10000,
        ymin=1, ymax=1000,
        legend pos=south east,
        grid style={line width=1pt, draw=black, opacity=0.1},
        xmajorgrids={true},
        ymajorgrids={true},
        log basis x=2
    ]
    \addplot
    [
        color=cp2,
        mark=square,
        thick
    ]
    coordinates {
    (365,20)(730,30.7)(1460,58.9)(2920,129.8)(8760,849)
    };
    \addplot
    [
        color=cp6,
        mark=square*,
        thick
    ]
    coordinates {
    (365,36.8)(730,28.6)(1460,30.4)(2920,32.48)(8760,57.1)
    };
    \addlegendentry{\PIPS}
    \addlegendentry{Hierarchical}
    \end{loglogaxis}
\end{tikzpicture}
\end{adjustbox}
\end{subfigure}%
\begin{subfigure}{.5\textwidth}
\caption{nb60\_to1}
\begin{adjustbox}{width=\textwidth}
\begin{tikzpicture}
    \begin{loglogaxis}
    [
        xlabel={\#blocks/\MPI processes},
        ylabel={Time to solution [s]},
        log ticks with fixed point,
        scale only axis,
        xmin=300, xmax=10000,
        ymin=1, ymax=1000,
        legend pos=south east,
        grid style={line width=1pt, draw=black, opacity=0.1},
        xmajorgrids={true},
        ymajorgrids={true},
        log basis x=2
    ]
    \addplot
    [
        color=cp2,
        mark=square,
        thick
    ]
    coordinates {
    (365,94.2)(730,88.7)(1460,148.4)(2920,270.4)(8760,849)
    };
    \addplot
    [
        color=cp6,
        mark=square*,
        thick
    ]
    coordinates {
    (365,130.4)(730,134)(1460,111.3)(2920,130.4)(8760,250.0)
    };
    \addlegendentry{\PIPS}
    \addlegendentry{Hierarchical}
    \end{loglogaxis}
\end{tikzpicture}
\end{adjustbox}
\end{subfigure}
\caption{Comparison of the original \PIPS and the hierarchical approach for different block sizes and constant amount of linking constraints.}
\label{fig:splitting}
\end{figure}

First, we see that the overall performance of the original approach for the given instances suffers with finer decomposition and only for the largest one can obtain some speedup when more blocks and processes are used.
The solution time here is dominated by the time spent in the Schur complement.
Only on the largest instance does the speed-up gained by using smaller diagonal blocks outweigh the growing Schur complement cost.

The behavior of \HSCA is similar only for the first instance.
On the other three instances there is some potential in using more and smaller diagonal blocks, a behavior directly coupled to \HSCA's sensitivity to the Schur complement size.
In \HSCA (in an ideal setting), the number of local linking constraints in each of the Schur complements grows only in $\bigO(\sqrt{n})$ (with $n$ denoting the number of local linking constraints), as compared to $\bigO(n)$ for the original approach.
\HSCA is thus able to scale ``further''.
Still, beyond a certain splitting, all instances experience some slow-down in \HSCA.

Comparing the runtime of both algorithms against each other, the original approach benefits from large diagonal blocks but experiences a significant slowdown with a growing number of linking constraints and shrinking diagonal blocks.
At the same time, \HSCA is worse at handling larger diagonal blocks and large Schur complements.
While it can outperform the original \PIPS for the smaller two instances, \HSCA falls behind as the diagonal block sizes increase.
For the hardest (in terms of diagonal block size) two instances, the original \PIPS cannot be outperformed by \HSCA.
In the other two instances, where the ratio of diagonal block size to the number of linking constraints favors \HSCA, the original algorithm consistently underperforms.
The optimal approach strongly depends on the individual layout of the given instance and the decomposition used.

There are two scenarios where we expect \HSCA, given sufficient compute resources, to consistently outperform the original \PIPS.
First, should a user not have access to a high-performance (commercial) sparse direct linear solver like \PARDISO but to (say) \MAFiftySeven, \HSCA should be picked over the original \PIPS.
A more lightweight solver will benefit from smaller matrices and scale well within \HSCA.
Second, in distributed but memory-bound environments, such as the \JUWELS supercomputer, large Schur complements can usually not be factorized practically or efficiently, and the only solution method available is the hierarchical extension.
In all other circumstances, the out-performance of either approach will be model-dependent.

\subsection*{Solution times of hierarchical approach on intractable models}

Lastly, we present experiments on models that the original implementation of \PIPS could not solve efficiently.
In \cref{tab:results}, we compare solution times obtained with \HSCA, the original \PIPS, and the three commercial solvers, \CPLEX 22.1.1.0, \GUROBI 11, and \COPT 7.2.3 each run with 32 threads.
Apart from the run time, we list model properties for each instance: the number of variables (\textit{vars.}), constraints (\textit{cons.}), nonzeros (\textit{nnz.}), and linking constraints (\textit{link. cons.}).
Variables, constraints, and nonzeros are given in millions (M).
We write MEM whenever \PIPS crashes due to memory limitations on the compute nodes.
For MISO\_DISP\_120 and MISO\_DISP\_488, we could not obtain results with the original \PIPS due to errors in underlying software packages.
For these same instances, Solver 2 could not obtain optimal solutions but instead returned suboptimal ones.
The respective times are marked with '*'.
The number of nodes used to obtain these results varies between 16 and 92; the column ``procs./nodes'' denotes for every instance the number of \MPI processes and compute nodes used.
For SIMPLE\_4, MISO\_DISP\_120, and MISO\_DISP\_488, we use more nodes than indicated by the number of \MPI processes ($\lceil \frac{\text{\#procs}}{48} \rceil$, 48 cores per node), as we additionally run with two OpenMP threads per process.
For SIMPLE\_4, \HSCA splits the Schur complement into 4 layers and 3 for the other instances.

\begin{table}[ht]
    \centering
    \sffamily
    \scriptsize
    \setlength{\tabcolsep}{2pt}
    \begin{tabular*}{\textwidth}{@{\extracolsep{\fill}}lrrrrrrrrrr@{}}
        \toprule
        &    \multicolumn{4}{c}{Size} &  & \multicolumn{5}{c}{Run time (seconds)}\\
        \cmidrule(lr){2-5}  \cmidrule(lr){7-11} 
        Instance         & vars.   & cons.  & nnz. & link. cons.    & procs./nodes &        \HSCA & \allcaps{PIPS} & Solver 1 & Solver 2 & Solver 3 \\
        \midrule
        SIMPLE\_1       &  5.79M &  5.26M &  20.94M &    525\,600 & 730/16 &       561 &     849 &     461 &     704 &     717 \\
        SIMPLE\_2       &  5.78M &  5.26M &  20.94M &    175\,200 & 730/16 &       121 &     278 &     479 &     710 &     744 \\
        SIMPLE\_3       & 59.80M & 51.60M & 205.57M &     71\,680 & 1\,024/43 &       922 &     400 &  4\,571 &  8\,076 &  5\,368 \\
        SIMPLE\_4       & 76.18M & 59.80M & 254.72M & 1\,679\,371 & 2\,190/48 &    2\,343 &     MEM &  6\,555 &  7\,902 &  6\,555 \\
        MISO\_EXP\_30   &  4.80M &  5.06M &  19.42M &     96\,388 & 2\,190/48 &        76 &     MEM &     795 &     855 &     647 \\
        MISO\_EXP\_120  & 14.56M & 15.35M &  57.90M &    311\,040 & 2\,190/48 &       741 &     MEM &  4\,690 & 12\,391 &  3\,235 \\
        MISO\_EXP\_240  & 22.37M & 26.40M &  98.39M & 1\,156\,398 & 1\,472/62 &    1\,392 &     MEM & 12\,252 & 25\,391 & 10\,624 \\
        MISO\_DISP\_120 & 11.32M & 12.11M &  45.99M &     48\,240 & 2\,190/92 &       48 &       - &     627 &  1\,037*&     397 \\
        MISO\_DISP\_488 & 25.44M & 11.79M & 120.09M &     54\,844 & 2\,190/92 &       206 &      - &  1\,373 &  3\,363*&  1\,060 \\
        \bottomrule
    \end{tabular*}%
    \caption{Computational results for large-scale instances (M=million).}
    \label{tab:results}
\end{table}

In line with the results presented in \cite{Rehfeldt2019_PIPSIPMpp}, \PIPS and \HSCA often outperform commercial software.
This is especially visible in real-world MISO instances and larger SIMPLE instances.
Overall, \HSCA can achieve speed-ups of up to a factor of 10 on the presented instances.
Similar to \cref{fig:splitting}, we see that for the instances SIMPLE\_2 and SIMPLE\_3, either approach can outperform the other.
This depends on the amount of linking constraints vs. the size of the diagonal blocks.
\PIPS cannot efficiently solve many instances as it often runs out of memory when factorizing the Schur complement.
None of the instances could be solved using the original \PIPS and \MAFiftySeven instead of \PARDISO.

We finally note that the performance of the commercial solvers often also strongly depends on their ability to presolve the instances.
They usually can reduce the problem size significantly more than \PIPS.
This has mainly two reasons.
First, \PIPS presolve is designed not to destroy the underlying problem structure of an instance \cite{Kempke19}.
This mostly leads to limitations when aggregating variables in between blocks.
Second, \PIPS current presolve only implements a handful of methods--far from what is available in modern commercial \LP solvers.
As a result, the presolved problems in \PIPS are often more challenging and exhibit higher redundancy than their commercial counterparts.
\section{Conclusions}
\label{sec:conclusions}

This paper presents \HSCA, a novel distributed approach for efficiently solving large arrowhead structured \LPs.
It builds on an efficient distributed factorization of the \KKT systems arising in \IPMs.
The approach exploits the primal–dual block-angular structure of our models to distribute and parallelize the underlying linear algebra.
Extending the original \PIPS implementation, we applied a multi-level Schur complement decomposition to deal with problems with up to 1.7 million linking constraints.
We showed the limits of our approach and discussed scenarios in which either implementation outperforms the other.
While \HSCA often outperforms commercial solvers and the original \PIPS, its performance is heavily problem-dependent.
In this respect, \HSCA complements the original \PIPS and should be considered an alternative when solving large-scale arrowhead \LPs.
Good scalability is the key feature that gives \PIPS and its hierarchical extension an edge over other commercial optimization software.
As demonstrated, most commercial software packages do not scale well with additional threads and cores.
As parallel hardware becomes increasingly affordable and efficient, scalable numerical algorithms are likely to outperform and eventually replace traditional solution methods, even when their total computational cost is higher.

While the first implementation of the hierarchical approach has already achieved promising results, several other research paths remain.
First, presolve is still limited, and a full presolve suite would improve both run time and stability of our solver and would make it more comparable to commercial software.
Second, we are exploring alternative linear solvers that may further improve \PIPS performance.
Third, as \HSCA is not limited to two-link constraints, a more elaborate structure detection could be further implemented.
In particular, one could account for linking variables and constraints linking $n$ consecutive blocks, leading to a more complicated Schur complement structure with additional ``global'' linking constraints in the subsystems.
The \IPM algorithm itself, discussed only briefly in this paper, also offers room for further improvement.
To this end, we plan on implementing a sequential version of \PIPS to robustify the \IPM and compare its performance against the performance of other codes.
Another promising direction is extending \HSCA to support approximate preconditioners instead of exact factorizations at each \IPM step.
Many of the preconditioners mentioned in \cite{BenziGolubLiesen_2005_NumericalSaddlePointSystems, DApuzzo2008_MutualImpactOfNumericalLAAndLargeScaleIPM, KarimSolomonik2022_EfficientPreconditionersForIPMViaSCStrategy} can directly profit from exploitation of the arrowhead structure and even be extended to a recursive arrowhead structure.

Last, many open questions remain to be addressed on the application side.
In the ongoing research project Peregrine\footnote{https://www.dlr.de/en/ve/research-and-transfer/projects/project-peregrine}, we are integrating automatic structure detection, similar to the work done in \cite{MartinDecomposingMatricesIntoBlocks1998}, into our software, as the modeler currently has to annotate a model’s block structure manually.
Lastly, we are preparing an open model library containing large-scale \ESMs from different contexts as well as distribution planning problems, which will also include the problems used in this paper.
We invite others to try out our software and contribute; the version used in this paper is available on \href{https://github.com/NCKempke/PIPS-IPMpp}{GitHub}.

\section*{Acknowledgements}
{\footnotesize The described research activities are funded by the Federal Ministry for Economic Affairs and Energy within the project UNSEEN (ID: 03EI1004C, 03EI1004D). The authors gratefully acknowledge the Gauss Centre for Supercomputing e.V. (www.gauss-centre.eu) for funding this project by providing computing time through the John von Neumann Institute for Computing (\allcaps{NIC}) on the \allcaps{GCS} Supercomputer \allcaps{JUWELS} at J{\"u}lich Supercomputing Centre (\allcaps{JSC}). We thank Charlie Vanaret for fruitful discussions and his fierce effort to make this paper more readable.}

\ifthenelse{\zibreport = 0}{
	\bibliographystyle{siamart_220329/siamplain}
}{
	\bibliographystyle{plain}
}
\bibliography{references}

\appendix

\section{Structure detection within the Schur complement}\label{apdx:sc_structure}

After realizing that many linking constraints only link two consecutive blocks as described in \cref{sec:linkcons_structure}, we can utilize this fact to detect structure in the Schur complement \cref{alg:factor_pips:SumSC} in \cref{alg:factor_pips} itself. Given the system matrix \cref{eq:linear_system_pips_perturbed_compressed} and defining
\begin{equation*}
    K_i^{-1} := \matr{ \tilde{K}^i_{1,1} & \tilde{K}^i_{1,2} \\ \tilde{K}^i_{2,1} & \tilde{K}^i_{2,2} }
\end{equation*}
the Schur complement contributions of each block $i$ have the following structure:
\begin{align*}
L_i^T K_i^{-1} L_i & = \matr{ 0 & A_i^T\\ 0 & 0 \\ \hat{\PIPSBl}_i & 0 \\ \PIPSBlglob_i & 0 }
    \matr{ \tilde{K}^i_{1,1} & \tilde{K}^i_{1,2} \\ \tilde{K}^i_{2,1} & \tilde{K}^i_{2,2} }
    \matr{ 0 & 0 & \hat{\PIPSBl}_i^T & \PIPSBlglob_i^T\\ A_i & 0 & 0 & 0 } \\
    & = \matr{ A_i^T \tilde{K}^i_{2,2} A_i & 0 & A_i^T \tilde{K}^i_{2,1} \hat{\PIPSBl}_i^T & A_i^T  \tilde{K}^i_{2,1} \PIPSBlglob_i^T \\
           0 & 0 & 0 & 0 \\
           \hat{\PIPSBl}_i \tilde{K}^i_{1,2} A_i & 0 & {\color{red}\hat{\PIPSBl}_i \tilde{K}^i_{1,1} \hat{\PIPSBl}_i^T} & \hat{\PIPSBl}_i
           \tilde{K}^i_{1,1} \PIPSBlglob_i^T  \\
           \PIPSBlglob_i \tilde{K}^i_{1,2} A_i  & 0 & \PIPSBlglob_i \tilde{K}^i_{1,1} \hat{\PIPSBl}_i^T & \PIPSBlglob_i \tilde{K}^i_{1,1} \PIPSBlglob_i^T }.
\end{align*}
The two-links imply additional structure on ${\color{red}\hat{\PIPSBl}_i \tilde{K}^i_{1,1} \hat{\PIPSBl}_i^T}$:
\begin{multline*}
    {\color{red}\hat{\PIPSBl}_i \tilde{K}^i_{1,1} \hat{\PIPSBl}_i^T} = \matr{ \vdots \\ 0 \\ \PIPSBlloc_{i-1}' \\ \PIPSBlloc_i \\ 0 \\ \vdots } \tilde{K}^i_{1,1}
    \matr{ \cdots & 0 & \PIPSBlloc_{i-1}'^T & \PIPSBlloc_i^T & 0 & \cdots }
    = \\
    \matr{\ddots & \vdots & \vdots & \vdots & \vdots & \iddots  \\
        \cdots   & 0 & 0 & 0 & 0 & \cdots \\
        \cdots   & 0 & \PIPSBlloc_{i-1}' \tilde{K}^i_{1,1} \PIPSBlloc_{i-1}'^T & \PIPSBlloc_{i-1}' \tilde{K}^i_{1,1} \PIPSBlloc_{i}^T & 0 & \cdots \\
        \cdots   & 0 & \PIPSBlloc_{i} \tilde{K}^i_{1,1} \PIPSBlloc_{i-1}'^T & \PIPSBlloc_{i} \tilde{K}^i_{1,1} \PIPSBlloc_{i}^T & 0 & \cdots \\
        \cdots   & 0 & 0 & 0 & 0 & \cdots \\
        \iddots  & \vdots & \vdots & \vdots & \vdots & \ddots  \\ }.
\end{multline*}
Summing up all individual Schur complement contributions (and after appropriate symmetric permutation), the final Schur complement is expressed as the following block matrix:
\begin{equation} \label{eq:schur_complement_blockstructure}
P^T  \big(	\sum_{i = 1}^{N} L_i^TK_i^{-1} L_i \big) P
= \matr{                 M_1  & H_1^T      &        &               & &    E_1^T & 0         \\
    H_1  &	M_2        & H_2^T &               & &    E_2^T  & 0         \\
    &	H_2        & \ddots &  \ddots      & &   \vdots  ~& \vdots   \\
    &            & \ddots &    &  H_{N - 2}^T &    ~~E_{N - 2}^T & 0           \\
    &         &   &   H_{N - 2}   & ~~M_{N-1}~~  &    ~~E_{N - 1}^T & 0 \\
    E_1  &  E_2 & \cdots &  E_{N-1} &  E_{N - 1}  &   M_0   ~ & 0 \\
    0 & 0 &   \cdots & 0 & 0 & 0 & 0
}.
\end{equation}
The submatrices in \cref{eq:schur_complement_blockstructure} are defined as
\begin{equation*}
     M_i := \PIPSBlloc_{i} \tilde{K}^i_{1,1} \PIPSBlloc_{i}^T + \PIPSBlloc_{i}' \tilde{K}^{i+1}_{1,1} \PIPSBlloc_{i}'^T,
\end{equation*}
\begin{equation}\label{eq:dense_border_part_sc}
    E_i :=  \matr{\PIPSBlglob_i  \tilde{K}^i_{1,1} \PIPSBlloc_i^T + \PIPSBlglob_{i+1} \tilde{K}^{i+1}_{1,1} \PIPSBlloc_i'^T  \\
    A_i^T \tilde{K}^i_{2,1} \PIPSBlloc_i^T + A_{i+1}^T \tilde{K}^{i+1}_{2,1} \PIPSBlloc_i'^T},
\end{equation}
for $i \in \{1,\dots,N-1\}$ and
\begin{equation*}
    H_i := \PIPSBlloc_{i+1} \tilde{K}^{i+1}_{1,1} \PIPSBlloc_{i}'^T,
\end{equation*}
\begin{equation}\label{eq:dense_corner_part_sc}
    M_0	:= \sum_{i=1}^{N} \matr{
    \PIPSBlglob_i  \tilde{K}^i_{1,1} \PIPSBlglob_i^T & \PIPSBlglob_i \tilde{K}^i_{1,2} A_i  \\
    A_i^T \tilde{K}^i_{2,1} \PIPSBlglob_i^T & A_i^T \tilde{K}^i_{2,2} A_i
    },
\end{equation}
for $i\in\{1,\dots,N-2\}$.
This structure is detected in \PIPS and used to precompute the sparsity pattern of the final Schur complement.
As long as no additional structure can be detected in the constraints, all blocks in \cref{eq:schur_complement_blockstructure} are dense; in particular, \cref{eq:dense_border_part_sc} and \cref{eq:dense_corner_part_sc} can be significantly large, depending on the number of linking constraints and global linking variables in the problem.
We reiterate an upper bound on the number of nonzeros in the Schur complement, as introduced in Observation 1 in \cite{Rehfeldt2019_PIPSIPMpp}:
\begin{equation*}
    \sum_{i=1}^{N-1} l_i^2 + 2 \sum_{i=1}^{N-2} l_i l_{i+1} + 2\sum_{i=1}^{N-1} l_i (m_\PIPSBlglob + n_0) + (m_\PIPSBlglob + n_0)^2.
\end{equation*}
Even though a two-link structure for linking variables could also be computed, linking variables are seldom local in the models considered. Therefore, this has not been implemented.

The Schur complements of systems without global linking constraints, as is the case for the inner linear systems defined by $K^{dense}$ in \cref{eq:linear_system_pips_inner_top} and ${K}^{inner}_i$ in \cref{eq:linear_system_pips_inner_lower} have the band diagonal form
\begin{equation} \label{eq:schur_complement_blockstructure_inner}
\matr{
\tilde{M}_1 & \tilde{H}_1^T &               &                 &                     \\
\tilde{H}_1 & \tilde{M}_2   & \tilde{H}_2^T &                 &                     \\
            & \tilde{H}_2   & \ddots        & \ddots          &                     \\
            &               & \ddots        &                 & \tilde{H}_{\tilde{N} - 2}^T \\
            &               &               & \tilde{H}_{\tilde{N}-2} & ~~\tilde{M}_{\tilde{N}-1}~~ \\
},
\end{equation}
with $\tilde{N}\in\N$ and $\tilde{M}_i \in \R^{\tilde{l}_i\times \tilde{l}_i}$ and $\tilde{H}_i \in \R^{\tilde{l}_{i+1}\times \tilde{l}_i}$ defined equivalently.
These Schur complements are mostly sparse.
The number of nonzeros of \cref{eq:schur_complement_blockstructure_inner} is thus bounded by $\sum_{i=1}^{\tilde{N}-1} \tilde{l}_i^2 + 2 \sum_{i=1}^{\tilde{N}-2} \tilde{l}_i \tilde{l}_{i+1}$.

\end{document}